\def\no{\if01}
\def\iftwelvept{\no}

\def\ifusepdf{\no}
\def\ifpsfont{\no}

\iftwelvept
\documentclass[leqno,12pt]{amsart}
\else
\documentclass[leqno,12pt]{amsart}
\fi
\usepackage{amssymb}
\usepackage{amscd}
\usepackage{latexsym}
\usepackage{verbatim}
\usepackage{mathrsfs}
\usepackage[all]{xy}

\ifusepdf
\usepackage{hyperref}
\else\fi
\ifpsfont
\usepackage[T1]{fontenc}
\usepackage{times}
\else\fi

\iftwelvept
\else\fi

\theoremstyle{plain}
\newtheorem{Theorem}{Theorem}[section]

\newtheorem{Proposition}[Theorem]{Proposition}
\newtheorem{Lemma}[Theorem]{Lemma}
\newtheorem{Corollary}[Theorem]{Corollary}

\theoremstyle{definition}

\newtheorem{Definition}[Theorem]{Definition}
\newtheorem{Remark}[Theorem]{Remark}
\newtheorem{Example}[Theorem]{Example}

\renewcommand{\theTheorem}{\arabic{section}.\arabic{Theorem}}
\renewcommand{\theClaim}{\arabic{section}.\arabic{Theorem}.\arabic{Claim}}
\renewcommand{\theequation}{\arabic{section}.\arabic{Theorem}.\arabic{Claim}}

\newcommand{\ZZ}{{\mathbb{Z}}}

\newcommand{\NNNN}{\operatorname{N}}

\newcommand{\I}{\infty}

\newcommand{\XXX}{{\mathscr{X}}}

\newcommand{\Sh}{\operatorname{Sh}}

\newcommand{\Rep}{\operatorname{Rep}}
\newcommand{\PRep}{\operatorname{PRep}}
\newcommand{\TTT}{{\mathcal{T}}}

\newcommand{\SM}{\operatorname{Sm}_S}

\newcommand{\MMM}{\mathcal{M}}
\newcommand{\NNN}{\mathcal{N}}

\newcommand{\Hom}{\operatorname{Hom}}

\newcommand{\Comp}{\operatorname{Comp}}
\newcommand{\Spec}{\operatorname{Spec}}

\newcommand{\Perf}{\operatorname{Perf}}
\newcommand{\Gal}{\operatorname{Gal}}

\newcommand{\Cor}{\operatorname{Cor}}
\newcommand{\SP}{\operatorname{Sp}}

\newcommand{\Mod}{\operatorname{Mod}}
\newcommand{\PMod}{\operatorname{PMod}}

\newcommand{\SSS}{\mathcal{S}}

\newcommand{\Map}{\operatorname{Map}}

\newcommand{\Fun}{\operatorname{Fun}}
\newcommand{\Alg}{\operatorname{Alg}}

\newcommand{\End}{\operatorname{End}}

\newcommand{\CoAlg}{\operatorname{CoAlg}}

\newcommand{\Aff}{\operatorname{Aff}}
\newcommand{\DM}{\mathsf{DM}}
\newcommand{\Grp}{\operatorname{Grp}}

\newcommand{\FIN}{\operatorname{Fin}_\ast}

\newcommand{\wCat}{\widehat{\textup{Cat}}_{\infty}}
\newcommand{\wsCat}{\widehat{\textup{Cat}}^{\textup{L},\textup{st}}_{\infty}}

\newcommand{\CAlg}{\operatorname{CAlg}}

\newcommand{\calC}{\mathcal{C}}

\newcommand{\LM}{\operatorname{LM}^\otimes}

\newcommand{\CHopf}{\operatorname{CHopf}}

\newcommand{\Ind}{\operatorname{Ind}}
\newcommand{\uCat}{\operatorname{Cat}_{\infty}}
\newcommand{\sCat}{\operatorname{Cat}^{\textup{st}}_{\infty}}

\newcommand{\Ass}{\operatorname{Ass}^{\otimes}}

\newcommand{\Aut}{\operatorname{Aut}}

\newcommand{\Sym}{\operatorname{Sym}}
\newcommand{\Proof}{{\sl Proof.}\quad}
\newcommand{\QED}{{\unskip\nobreak\hfil\penalty50\quad\null\nobreak\hfil
{$\Box$}\parfillskip0pt\finalhyphendemerits0\par\medskip}}

\begin{document}

\title{Tannakization in derive algebraic geometry}

\author{Isamu Iwanari}

\thanks{The author is partially supported by Grant-in-Aid for Scientific Research,
Japan Society for the Promotion of Science.}

\begin{abstract}
We give a universal construction of a derived affine group scheme
and its representation category
from a symmetric monoidal $\infty$-category, which we shall call the tannnakization of
a symmetric monoidal $\infty$-category.
This can be viewed as an $\infty$-categorical generalization of
 the work of Joyal-Street \cite{JS} and Nori.
We then apply it to the stable $\infty$-category of mixed motives
equipped with the realization functor of a mixed Weil cohomology
and obtain a derived motivic Galois group whose representation
category has a universality, and which represents the automorphism
group of the realization functor.
Also, we present basic properties
of derived affine group schemes in Appendix.
\end{abstract}

\address{Mathematical Institute, Tohoku University, Sendai, Miyagi, 980-8578,
Japan}

\email{iwanari@math.tohoku.ac.jp}

\maketitle

\section{Introduction}

Grothendieck has developed the theory of Galois categories \cite{SAG1},
and Saavedra and Deligne-Milne have studied the theory of tannakian categories
\cite{Sa}, \cite{DM} which generalizes the classical
Tannaka duality \cite{Ta}
by the categorical and algebro-geometric method.
These are beautiful duality theories in their own right on one hand,
one of important aspects of these theories is the
role as the powerful machine by which we can derive invariants from 
abstract categories on the other hand.
For example, the \'etale fundamental groups
of schemes and Picard-Vessiot Galois groups
were constructed by means of these theories.
Joyal-Street \cite{JS} and Nori gave
the machinery which approximates symmetric monoidal categories
and graphs
with (neutral) tannakian abelian categories (the braided case was also treated
in \cite{JS}).
This machinery is powerful: Joyal-Street applied it to quantum groups, and Nori used it to construct the Nori's category of motives (see e.g.
\cite{Ara}).
We here informally call this approximation the tannakization of categories.

The first main purpose is to construct
tannakization in the setting of higher categories, i.e.
$\infty$-categories.
In this introduction, by an $\infty$-category
we informally mean a (weak) higher category, in which all $n$-morphisms are weakly
invertible for $n>1$ (cf. \cite{Be}). (There
are several theories which provide ``models'' of such categories.
We use the machinery of quasi-categories from the next
Section.)
Let $\mathcal{C}^{\otimes}$ be a symmetric monoidal small
$\infty$-category.
For a commutative ring spectrum $R$, we let $\Mod_R^\otimes$
be the symmetric monoidal $\infty$-category of $R$-module spectra.
Let $\PMod_R^\otimes$ be the symmetric monoidal full subcategory
of $\Mod_R^\otimes$ spanned by dualizable objects (cf. Section 2).
Let $\CAlg_R$ be the $\infty$-category of commutative $R$-ring spectra.
Let $\omega:\mathcal{C}^\otimes \to \PMod_R^\otimes$ be a symmetric monoidal
functor.
Then our result can be roughly stated as follows (see Theorem~\ref{supermain}):

\begin{Theorem}
\label{intro1}
There are a derived affine
group scheme $G$ over $R$ (explained below)
and a symmetric monoidal functor
$u:\mathcal{C}^\otimes\to \PRep_G^\otimes$ which makes the outer triangle in
\[
\xymatrix{
     &  \PMod_G^\otimes \ar@{..>}[d] \ar[ddr]^{\textup{forget}}  &     \\
   & \PRep_H^\otimes \ar[rd]_{\textup{forget}} &  \\
\mathcal{C}^\otimes \ar[ru] \ar[ruu]^{u} \ar[rr]_\omega & & \PMod_R^\otimes
}
\]
commute in the $\infty$-category of symmetric monoidal $\infty$-categories
(here $\PRep_G^\otimes$ is the symmetric monoidal
$\infty$-category of dualizable
$R$-module spectra equipped with $G$-actions)
such that these possess the following universality:
for any inner triangle consisting of solid arrows in the above diagram
where $H$ is a derived affine group scheme,
there exists a unique (in an appropriate sense)
morphism $f:H\to G$ of derived affine group
schemes which induces $\PRep_G^\otimes\to \PRep_H^\otimes$
(indicated by the dotted arrow) filling the above diagram.
Moreover, the automorphism group of $\omega$ is represented by $G$.
\end{Theorem}

For simplicity, we usually
refer to the pair $(G, u:\mathcal{C}^\otimes\to \PRep_G^\otimes)$ as the tannakization.
By Theorem~\ref{intro1} we can obtain ``Tannaka-Galois type
invariants''
in the quite general setting.
A derived group scheme is an analogue of
group schemes in derived algebraic geometry.
This notion plays an important role in this paper.
To understand why this notion comes in,
let us recall that stable $\infty$-categories
are enriched over spectra
(cf. \cite[Section 2]{BGT}).
It leads us to
consider commutative Hopf ring spectra which are
the spectra version of commutative Hopf algebra.
Put another way, from an intuitive point of view, pro-algebraic groups
(i.e. affine group schemes) appears in the formulation of
classical Tannaka duality
since the automorphisms of finite-dimensional vector
spaces are representable by algebraic groups.
Similarly, the automorphisms of compact spectra (or a bounded
complexes of finite dimensional vector spaces) are representable by derived affine group schemes.
The fundamental and comprehensive works on
derived algebraic geometry by To\"en-Vezzosi
\cite{HAG2}, Lurie \cite{DAGn} provide a natural home in which
one can realize this idea.
For example, the functor $\omega$ can possess higher automorphisms.
The derived affine group scheme $G$
captures all these higher data.

We would like to stress that we impose only weak natural conditions
on $\mathcal{C}^\otimes$ and $\omega$ in Theorem~\ref{intro1}.
Consequently, it is applicable also to situations in which
$\mathcal{C}^\otimes$
 seems ``non-tannakian''.
Typical examples are
$\mathcal{C}^\otimes=\PMod_A^\otimes$
with $A$ arbitrary.
Even in the case, our tannakization provides meaningful invariants.
In a separate paper \cite{Bar}, we prove that our tannakization
includes bar construction of an augmented
commutative ring spectrum and its equivariant versions
as a special case.
Therefore our tannakization can be also viewed as a generalization of
bar constructions and equivariant bar constructions.

\vspace{2mm}

Our motivation comes from various important and interesting examples
which live in the realm of $\infty$-categories.
For example, the triangulated
category of mixed motives, due to Hanamura, Levine and Voevodsky,
is of great interest in the view of a tannakian theory
for higher categories.
The
category of mixed motives
has a natural formulation of symmetric monoidal
stable $\infty$-category.
The stable $\infty$-category is equipped with realization functors of mixed
Weil
cohomology theories.
One of important examples of stable $\infty$-categories
which recently appeared might be a symmetric monoidal stable $\infty$-category of noncommutative
motives by Blumberg-Gepner-Tabuada \cite{BGT}, that is 
the natural and universal domain for localizing (or additive) invariants such as
algebraic K-theory, topological Hochschild homology and topological cyclic homology.
As for an example which is not ``algebraic'' one,
the stable $\infty$-category
of the perfect complexes on a topological space gives us an
tannakian invariant.

\vspace{2mm}

From Section 5, we then switch to applications
to examples.
We will construct a derived motivic Galois group for mixed motives.
We note that
for our construction we do not need a conjectural
motivic $t$-structure
(see also Remark~\ref{abelianmotives} on this point).
Let $\mathbf{K}$ be a field of characteristic zero
and let $H\mathbf{K}$ denote the Eilenberg-MacLane spectrum.
Let $\mathsf{DM}^\otimes:=\mathsf{DM}^\otimes(k)$ be the $H\mathbf{K}$-linear symmetric monoidal stable $\infty$-category
of mixed motives over a perfect field $k$
(see Section 5).
Let $\mathsf{DM}^\otimes_\vee$ be the symmetric monoidal
full subcategory spanned by dualizable objects in $\mathsf{DM}^\otimes$.
In $\mathsf{DM}^\otimes$, dualizable objects coincide with compact objects.
The homotopy category of $\mathsf{DM}_\vee$ can be identified with
the $\mathbf{K}$-linear
triangulated category of geometric
motives $DM_{gm}(k)$ constructed by
Voevodsky (see e.g. \cite{MVW}, \cite{Vtri}), which is anti-equivalent to Hanamura's category \cite{HaM}
and
Levine's category \cite{LeM} (with rational coefficients).
Let $E$ be a mixed Weil (cohomology) theory
with coefficients $\mathbf{K}$
in the sense of \cite{CD2}.
For example, $l$-adic \'etale cohomology, Betti cohomology,
de Rham cohomology
and rigid cohomology
give mixed Weil theories.
Then we can construct the homological realization functor
\[
\mathsf{R}_E:\mathsf{DM}^\otimes_\vee\longrightarrow\PMod_{H\mathbf{K}}^\otimes,
\]
that is
a symmetric monoidal exact functor (see Section 5).
Note that the homotopy category of
$\PMod_{H\mathbf{K}}^\otimes$
can be regarded as the triangulated category of
bounded complexes of $\mathbf{K}$-vector speces with
finite dimensional cohomology groups.
Applying Theorem~\ref{intro1} to the realization functor of a mixed Weil
cohomology theory
we obtain (cf. Definition~\ref{mgg}, Theorem~\ref{derivedgroup}):

\begin{Theorem}
\label{intro2}
The realization functor $\mathsf{R}_E:\mathsf{DM}^\otimes_\vee\to \PMod_{H\mathbf{K}}^\otimes$
gives rise to the tannakization
$(\mathsf{MG}_E=\Spec B_E,\mathsf{DM}^\otimes_\vee\to \PRep_{\mathsf{MG}_E}^\otimes)$
over $H\mathbf{K}$ described in Theorem~\ref{intro1}.
Here $B_E$ is a commutative differential graded $\mathbf{K}$-algebra.
\end{Theorem}

By the universality and representability,
we shall propose $\mathsf{MG}_E$ as a
(derived) motivic Galois group of mixed
motives.
By a truncation procedure we can also
extract the underived motivic Galois group
$MG_E$, which is an ordinary affine group scheme,
from $\mathsf{MG}_E$ (cf. Theorem~\ref{ordinarygroup}).
Apart from the universality and representability of $\mathsf{MG}_E$,
our derived group scheme $\mathsf{MG}_E$
can be thought of as a natural generalization of so-called
motivic Galois group $MTG$ for mixed Tate motives constructed by
Bloch-Kriz, Kriz-May, Levine \cite{BK}, \cite{KM}, \cite{Lev}.
To explain this, we would like to invite the reader's attention
to the results obtained in \cite{Bar}.
Let $\mathsf{DTM}_\vee\subset \mathsf{DM}_\vee$ be the stable
$\infty$-category of mixed Tate motives, that is,
the stable idempotent complete subcategory generated by Tate
objects $\{\mathbf{K}(n)\}_{n\in \ZZ}$ (see \cite[Section 6]{Bar} for more
details).
The full subcategory $\mathsf{DTM}_\vee$ forms a symmetric monoidal
$\infty$-category $\mathsf{DTM}^\otimes$, whose symmetric monoidal structure is
induced by that of $\mathsf{DM}^\otimes$.
In \cite{Bar}, we prove the comparison results which can be informally
summarized as
follows:

\begin{Theorem}[\cite{Bar}]
\label{intro3}
\begin{enumerate}
\renewcommand{\labelenumi}{(\roman{enumi})}
\item Let $\mathsf{MTG}$ be a derived affine group
scheme over $H\mathbf{K}$ obtained as the tannakization
of the composite $\mathsf{R}_T:\mathsf{DTM}_\vee^\otimes\hookrightarrow \mathsf{DM}^\otimes\stackrel{\mathsf{R}_E}{\to} \PMod_{H\mathbf{K}}^\otimes$
(we omit the subscript $E$).
Then 
$\mathsf{MTG}$ is equivalent to
a derived affine group scheme obtained from
the
$\mathbb{G}_m$-equivariant bar construction of a commutative
differential graded $\mathbf{K}$-algebra $\overline{Q}$ equipped with
$\mathbb{G}_m$-action.
That is to say, it is the \v{C}ech nerve of a morphism of derived
stacks $\Spec H\mathbf{K}\to [\Spec \overline{Q}/\mathbb{G}_m]$
(cf. Appendix Example~\ref{barcon} and \cite{Bar}). The complex $\overline{Q}$
is described in term of Bloch's cycle complexes.

\item Suppose that Beilinson-Soul\'e vanishing conjecture holds for $k$.
Let $MTG$ be the Tannaka
dual of the tannakian category of mixed Tate motives; the heart of motivic $t$-structure on $\mathsf{DTM}_\vee$ (constructed under the vanishing conjecture,
see \cite{Lev}, \cite{KM}, \cite[Section 7]{Bar}).
Then the affine group scheme $MTG$ is the underlying group scheme
(cf. Appendix A.4 or \cite[7.3]{Bar}) of $\mathsf{MTG}$.

\item Let $\mathsf{Art}^\otimes$ be the symmetric monoidal
stable idempotent complete
full subcategory generated by motives of
smooth zero-dimensional varieties, i.e. Artin motives.
Then the tannakization of $\mathsf{Art}^\otimes$ equipped with
a realization functor is the absolute Galois group $\Gal(\bar{k}/k)$.
\end{enumerate}
\end{Theorem}

This result links the works on mixed Tate motives
in \cite{BK}, \cite{KM}, \cite{Lev} and the classical Galois theory
to our results.
Adams graded bar constructions (that is, $\mathbb{G}_m$-equivariant bar constructions) are the fundamental tools in \cite{BK} and \cite{KM},
and
the central theme of \cite{Bar} is to compare bar constructions
and tannakizations.
In a sense,
the aspect of tannakiziation, that is Theorem~\ref{intro1},
as a generalization of bar constructions
allows us to construct a motivic Galois group of all mixed
motives.
In addition, it is worth mentioning that
Theorem~\ref{intro1} can be applied to any symmetric monoidal
full subcategory in $\mathsf{DM}^\otimes$.

We would like to emphasize
that higher category theory ($\infty$-categories) and
derived algebraic geometry provide a natural and nice
framework for our purposes.
For a commutative ring spectrum $A$,
the homotopy category of $\PMod_A^\otimes$ (or $\Mod_A^\otimes$)
forms a triangulated category equipped with a symmetric
monoidal structure.
However, if we work with triangulated categories (to prove Theorem~\ref{intro1} in particular, representability),
we encounter several technical problems including
the problem concerning
the absence of descent of morphisms in the homotopy category of $\PMod_A$.
It turns out
that $\infty$-categories give us an appropriate theory.
Also, we should like to refer the reader to
the recent preprints
\cite[VIII]{DAGn} \cite{Wall} and our previous work
\cite{FI} building on tannakian philosophy
in higher category theory.

The notion of derived (affine) group schemes is placed at the important part
of
our work.
We hereby decide to give the basic theory of derived affine group schemes
in Appendix. We also refer the reader to \cite{TH} and \cite{Spi}
for other accounts of related notions.

\vspace{2mm}

This paper is roughly organized as follows.
In Section 2, we fix notation and convention.
In Section 3 we give preliminaries which we need
Section 4.
Section 4 is devoted to the proof of
Theorem~\ref{intro1}.
Section 5 contains the construction of our motivic Galois group;
we construct the realization functor in the setting
of $\infty$-categories and apply Theorem~\ref{intro1}
to obtain a derived motivic Galois group
associated to the stable $\infty$-category of mixed motives.
In Section 6, we present some other examples without proceeding into
detail.
One example given in Section 6 is
the $\infty$-category of perfect complexes on a topological space $S$.
With rational coefficients, we expect that the associated derived affine group
is closely related to the rational homotopy theory.
It would yield a conceptual understanding of the rational homotopy theory
as an example of the tannakian philosophy.
In Appendix we present basic definitions and results
concerning derived
group schemes.

\section{Notation and Convention}

We fix notation and convention.

{\it $\infty$-categories}.
In this paper, we use theory of quasi-categories.
A quasi-category is a simplicial set which
satisfies the weak Kan condition of Boardman-Vogt:
A quasi-category $S$ is a  simplicial set
such that for any $0< i< n$ and any diagram
\[
\xymatrix{
\Lambda^{n}_i \ar[r] \ar[d] & S \\
\Delta^n \ar@{..>}[ru] & 
}
\]
of solid arrows, there exists a dotted arrow filling the diagram.
Here $\Lambda^{n}_i$ is the $i$-th horn and $\Delta^n$
is the standard $n$-simplex.
The theory of quasi-categories from higher categorical viewpoint
has been extensively developed by Joyal and Lurie.
Following \cite{HTT} we shall refer to quasi-categories
as $\infty$-categories.
Our main references are \cite{HTT}
 and \cite{HA}
(see also \cite{Jo}, \cite{DAGn}).
We often refer to a map $S\to T$ of $\infty$-categories
as a functor. We call a vertex in an $\infty$-category $S$
(resp. an edge) an object (resp. a morphism).
For the rapid introduction to $\infty$-categories, we
refer to \cite[Chapter 1]{HTT}, \cite[Section 2]{FI}.
It should be emphasized that
there are several alternative theories such as Segal categories,
complete Segal spaces, simplicial categories, relative categories,...
etc.
For the quick survey on various approaches to $(\infty,1)$-categories
and their relations,
we refer the reader to \cite{Be}.

\begin{itemize}

\item $\Delta$: the category of linearly ordered finite sets (consisting of $[0], [1], \ldots, [n]=\{0,\ldots,n\}, \ldots$)

\item $\Delta^n$: the standard $n$-simplex

\item $\textup{N}$: the simplicial nerve functor (cf. \cite[1.1.5]{HTT})

\item $\mathcal{C}^{op}$: the opposite $\infty$-category of an $\infty$-category $\mathcal{C}$

\item Let $\mathcal{C}$ be an $\infty$-category and suppose that
we are given an object $c$. Then $\mathcal{C}_{c/}$ and $\mathcal{C}_{/c}$
denote the undercategory and overcategory respectively (cf. \cite[1.2.9]{HTT}).

\item $\operatorname{Cat}_\infty$: the $\infty$-category of small $\infty$-categories in a fixed Grothendieck universe (cf. \cite[3.0.0.1]{HTT})

\item $\wCat$: $\infty$-category of $\infty$-categories

\item $\SSS$: $\infty$-category of small spaces (cf. \cite[1.2.16]{HTT})

\item $\textup{h}(\mathcal{C})$: homotopy category of an $\infty$-category (cf. \cite[1.2.3.1]{HTT})

\item $\Fun(A,B)$: the function complex for simplicial sets $A$ and $B$

\item $\Fun_C(A,B)$: the simplicial subset of $\Fun(A,B)$ classifying
maps which are compatible with
given projections $A\to C$ and $B\to C$.

\item $\Map(A,B)$: the largest Kan complex of $\Fun(A,B)$ when $A$ and $B$ are $\infty$-categories,

\item $\Map_C(A,B)$: the simplicial subset of $\Map(A,B)$ classifying
maps which are compatible with
given projections $A\to C$ and $B\to C$.

\item $\Map_{\mathcal{C}}(C,C')$: the mapping space from an object $C\in\mathcal{C}$ to $C'\in \mathcal{C}$ where $\mathcal{C}$ is an $\infty$-category.
We usually view it as an object in $\mathcal{S}$ (cf. \cite[1.2.2]{HTT}).

\end{itemize}

{\it Symmetric monoidal $\infty$-categories and spectra}.
We employ the theory of symmetric monoidal $\infty$-categories
developed in \cite{HA}.
We refer to \cite{HA} for its generalities.
Let $\FIN$ be the category of marked finite sets
(our notation is slightly different from \cite{HA}).
Namely, objects are marked finite sets and
a morphism from $\langle n\rangle_*:=\{1<\cdots <n\}\sqcup\{*\}\to \langle m\rangle_*:=\{1<\cdots <m\}\sqcup\{*\}$
is a (not necessarily order-preserving) map of finite sets which preserves the distinguished points $*$.
Let $\alpha^{i,n}:\langle n\rangle_*\to \langle 1\rangle_*$ be
a map such that $\alpha^{i,n}(i)=1$ and $\alpha^{i,n}(j)=\ast$ if $i\neq j\in \langle n\rangle_*$.
A symmetric monoidal category is a coCartesian fibration (cf. \cite[2.4]{HTT})
$p:\mathcal{M}^\otimes\to \NNNN(\FIN)$ such that
for any $n\ge0$, $\alpha^{1,n}\ldots \alpha^{n,n}$ induce
an equivalence $\mathcal{M}^\otimes_{n}\to (\mathcal{M}^\otimes_1)^{\times n}$
where $\mathcal{M}^\otimes_{n}$ and $\mathcal{M}^\otimes_{1}$
are fibers of $p$ over $\langle n\rangle_*$ and $\langle 1\rangle_*$
respectively.
A symmetric monoidal functor is a map $\mathcal{M}^\otimes\to \mathcal{M}^{'\otimes}$ of coCartesian fibrations
over $\NNNN(\FIN)$,
which carries coCartesian edges to coCartesian edges.
Let $\uCat^{\Delta,\textup{sMon}}$
be the simplicial category of symmetric monoidal $\infty$-categories
in which morphisms are symmetric monoidal functors.
Hom simplicial sets are given by those defined in \cite[3.1.4.4]{HTT}.
Let $\uCat^{\textup{sMon}}$ 
be the simplicial nerve of $\uCat^{\Delta,\textup{sMon}}$ (see \cite[2,1.4.13]{HA}).

There are several approaches to a ``good''theory of commutative ring spectra.
Among these,
we employ the theory of spectra and commutative ring spectra
developed in \cite{HA}.

We list some of notation.

\begin{itemize}

\item $\mathbb{S}$: the sphere spectrum



\item $\Mod_A$: $\infty$-category of $A$-module spectra
for a commutative ring spectrum $A$

\item $\PMod_A$: the full subcategory of $\Mod_A$ spanned by compact objects
(in $\Mod_A$, an object is compact if and only if it is dualizable, see \cite{BFN}) .
We refer to objects in $\PMod_A$ as perfect $A$-module (spectra).

\item Let $\mathcal{M}^\otimes\to \mathcal{O}^\otimes$ be a fibration of $\infty$-operads. We denote by
$\Alg_{/\mathcal{O}^\otimes}(\mathcal{M}^\otimes)$ the $\infty$-category of algebra objects (cf. \cite[2.1.3.1]{HA}).  We often write $\Alg(\MMM^\otimes )$ or 
$\Alg(\MMM)$ for $\Alg_{/\mathcal{O}^\otimes}(\MMM^\otimes)$.
Suppose that $\mathcal{P}^\otimes\to \mathcal{O}^\otimes$ is a map of $\infty$-operads. $\Alg_{\mathcal{P}^\otimes/\mathcal{O}^\otimes}(\mathcal{M}^\otimes)$: $\infty$-category of $\mathcal{P}$-algebra objects.

\item $\CAlg(\mathcal{M}^\otimes)$: $\infty$-category of commutative
algebra objects in a symmetric
monoidal $\infty$-category $\mathcal{M}^\otimes\to \NNNN(\FIN)$.

\item $\CAlg_R$: $\infty$-category of commutative
algebra objects in the symmetric monoidal $\infty$-category $\Mod_R^\otimes$
where $R$ is a commutative ring spectrum. When $R=\mathbb{S}$, we set
$\CAlg=\CAlg_{\mathbb{S}}$. The $\infty$-category $\CAlg_R$ is equivalent to
the undercategory $\CAlg_{R/}$ as an $\infty$-category.

\item $\Mod_A^\otimes(\mathcal{M}^\otimes)\to \NNNN(\FIN)$: symmetric monoidal
$\infty$-category of
$A$-module objects,
where $\mathcal{M}^\otimes$
is a symmetric monoidal $\infty$-category such that (1)
the underlying $\infty$-category admits a colimit for any simplicial diagram, and (2)
its tensor product functor $\mathcal{M}\times\mathcal{M}\to \mathcal{M}$
preserves
colimits of simplicial diagrams separately in each variable.
Here $A$ belongs to $\CAlg(\mathcal{M}^\otimes)$
(cf. \cite[3.3.3, 4.4.2]{HA}).

\end{itemize}

 Let $\mathcal{C}^\otimes$ be the symmetric monoidal $\infty$-category.
We usually denote, dropping the subscript $\otimes$,
by $\mathcal{C}$ its underlying $\infty$-category.
We say that an object $X$ in $\mathcal{C}$
is dualizable if there exist an object
$X^\vee$ and two morphisms $e:X\otimes X^\vee\to \mathsf{1}$ and
$c:\mathsf{1} \to X\otimes X^\vee$
with $\mathsf{1}$ a unit such that the composition
\[
X \stackrel{\textup{Id}_X\otimes c}{\longrightarrow} X\otimes X^\vee\otimes X \stackrel{e\otimes\textup{Id}_X}{\longrightarrow} X
\]
is equivalent to the identity, and 
\[
X^\vee \stackrel{c\otimes\textup{Id}_{X^\vee}}{\longrightarrow} X^\vee \otimes X\otimes X^\vee \stackrel{\textup{Id}_{X^\vee}\otimes e}{\longrightarrow} X^\vee
\]
is equivalent to the identity.
The symmetric monoidal structure of $\mathcal{C}$ induces
that of the homotopy category
$\textup{h}(\mathcal{C})$.
If we consider $X$ to be an object also
in $\textup{h}(\mathcal{C})$,
then $X$ is dualizable in $\mathcal{C}$ if and only if $X$
is dualizable in $\textup{h}(\mathcal{C})$.
For example, for $R\in \CAlg$, compact and dualizable objects
coincide in the symmetric monoidal $\infty$-category
$\Mod_R^\otimes$ (cf. \cite{BFN}).

\vspace{1mm}

\section{Basic definitions and geometric systems}

In this Section, we prepare some notions which we need in the next Section.

The $\infty$-category $\uCat$ of small $\infty$-categories
has the
symmetric monoidal structure determined by the Cartesian product
$\mathcal{C}\times \mathcal{D}$.
We denote by $\CAlg(\uCat)$ the $\infty$-category of commutative algebra
(monoid)
objects in the symmetric monoidal $\infty$-category $\uCat$.
A symmetric monoidal $\infty$-category
can be identified with a commutative algebra (monoid) object
in $\uCat$; there is a natural categorical equivalence
$\uCat^{\textup{sMon}}\simeq\CAlg(\uCat)$.
If $\mathcal{A}^\otimes, \mathcal{B}^\otimes\in \CAlg(\uCat)$,
we write $\Map^\otimes(\mathcal{A}^\otimes,\mathcal{B}^\otimes)$
for $\Map_{\CAlg(\uCat)}(\mathcal{A}^\otimes,\mathcal{B}^\otimes)$.

\vspace{1mm}

{\it Geometric $R$-system.}
We introduce the notion of geometric $R$-systems.

\begin{Definition}
Let $\mathcal{T}^\otimes:\CAlg_R\to \uCat^{\textup{sMon}}\simeq \CAlg(\uCat)$ be a functor satisfying the following
properties:
\begin{itemize}
\item[(A1)] Let $\mathcal{T}:\CAlg_R\to \CAlg(\uCat)\to \uCat$
be the composition with the forgetful functor.
For any $A$,
$\mathcal{T}(A)$ is stable and $\mathcal{T}(A)\to \mathcal{T}(B)$ is exact
for any $A\to B$.
For any $T\in \mathcal{T}(R)$, the automorphism group functor
$\Aut(T):\CAlg_R\to \mathcal{S}$, which will be defined below, is representable by a derived affine scheme over $R$.

\item[(A2)] For any $T, T'\in \mathcal{T}(R)$, the hom functor
$\Hom(T,T'):\CAlg_R\to \mathcal{S}$, which will be defined below, is representable by a derived affine scheme over $R$.

\end{itemize}
If (A1) and (A2) hold, we
refer to $\mathcal{T}^\otimes$ as a geometric $R$-system.
\end{Definition}
We here define $\Hom(T,T'):\CAlg_R\to \mathcal{S}$ as follows.
Let $\theta_{\Delta^1},\theta_{\partial\Delta^1},\theta_{\phi}:\uCat\rightarrow \mathcal{S}$ be the functors corresponding to $\Delta^1$, $\partial\Delta^1$ and the empty category $\phi$ respectively
via the Yoneda embedding $\uCat^{op}\subset \Fun(\uCat,\mathcal{S})$.
The inclusion $\partial\Delta^1\hookrightarrow \Delta^1$
induces $\theta_{\Delta^1}\to \theta_{\partial\Delta^1}$.
Note that $\theta_{\phi}$ is equivalent to the
constant functor whose value is the contractible space.
The functor $\theta_{\partial\Delta^1}$
is equivalent to the 2-fold product of the functor
$\uCat\to \mathcal{S}$ which carries an $\infty$-category $\mathcal{A}$
to the largest Kan complex $\mathcal{A}^{\simeq}$
(this functor can be constructed as the functor corepresentable by
$\Delta^0$).
Therefore, if we let $\mathcal{F}\to \CAlg_R$
be a left fibration
corresponding to $\CAlg_R\stackrel{T}{\to} \uCat\stackrel{\theta_{\Delta^0}}{\to} \mathcal{S}$, then
giving $\theta_{\phi}\to \theta_{\partial\Delta^1}$
amounts to giving two sections of $\mathcal{F}\to \CAlg_R$.
In order to construct $\theta_{\phi}\to \theta_{\partial\Delta^1}$
from $T$ and $T'$, we give (ordered) two sections $\CAlg_R \to \mathcal{F}$.
By \cite[3.3.3.4]{HTT}, a section corresponds to
an object in the limit $\lim\mathcal{T}(A)$
of $\mathcal{T}:\CAlg_R\to \uCat$.
Hence the images of $T$ and $T'$ in $\lim\mathcal{T}(A)$
give rise to $\theta_{\phi}\to \theta_{\partial\Delta^1}$.
We define $\Hom(T,T')$ to be the fiber product
$\theta_{\phi}\times_{\theta_{\partial\Delta^1}}\theta_{\Delta^1}$
in $\Fun(\CAlg_R,\mathcal{S})$.
For any $A\in \CAlg_R$, $\Hom(T,T')(A)$ is equivalent to (homotopy)
fiber product
\[
\{(T\otimes_RA,T'\otimes_RA)\}\times_{\Map(\partial\Delta^1, \mathcal{T}(A))}\Map(\Delta^1, \mathcal{T}(A))
\]
in $\mathcal{S}$,
where $T\otimes_RA$ and $T'\otimes_RA$ denote the images of $T$ and $T'$
in $\mathcal{T}(A)$ respectively. It is the mapping space from $T\otimes_RA$
to $T'\otimes_RA$.
If $T=T'$, we write $\End(T)$ for $\Hom(T,T)$.
We let $\Aut(T)$ be the functor $\CAlg_R\to \mathcal{S}$
obtained by restricting objects in $\End(T)(A)$ to automorphisms
for each $A$ (one can do this procedure by using
corresponding left fibration).

The followings are examples of geometric $R$-systems.
In
the next Section we will prove that these examples are geometric $R$-systems.

\begin{Example}
\label{example1}
Let $\Theta:\CAlg_R\to \CAlg(\uCat)$ be the functor which
carries $A$ to $\PMod_A^\otimes$ and carries $A\to B$ to
the base change functor $\PMod_A^\otimes\to \PMod_B^\otimes$.
We can obtain this functor $\Theta$ as follows.
By virtue of \cite[6.3.5.18]{HA}, we have
$\CAlg_R\to \CAlg(\widehat{\textup{Cat}}_\infty)_{\Mod_R^\otimes/}$
which carries $A$ to $\Mod_A^\otimes$.
Composing with the forgetful functor
$\CAlg(\widehat{\textup{Cat}}_\infty)_{\Mod_R^\otimes/}\to \CAlg(\widehat{\textup{Cat}}_\infty)$ and restricting $\Mod_A^\otimes$ to $\PMod_A^\otimes$
we have $\Theta:\CAlg_R\to \CAlg(\uCat)$.
This is a geometric $R$-system.
\end{Example}

\begin{Example}
\label{example2}
Let $S$ be a (small)
Kan complex.
Let $f:\uCat\to \uCat$ be the colimit-preserving
functor which is determined by
$(-)\times S$.
Namely, $f$ carries $\mathcal{C}$ to $\mathcal{C}\times S$.
Its right adjoint functor $g:\uCat\to \uCat$
carries $\mathcal{C}$ to $\Fun(S,\mathcal{C})$.
(To obtain this adjoint, consider the adjunction $(-)\times S:\operatorname{Set}_\Delta\rightleftarrows \operatorname{Set}_\Delta:\Fun(S,-)$, where
$\operatorname{Set}_\Delta$ denotes the category of simplicial sets.
If both $\operatorname{Set}_\Delta$ are endowed with Joyal model
structure \cite[2.2.5.1]{HTT}, then this adjunction is a Quillen adjunction
by \cite[2.2.5.4]{HTT}. It gives rise to
the required adjunction)
Then $\Fun(\NNNN(\FIN),\mathcal{C})\to \Fun(\NNNN(\FIN),\mathcal{C})$
induced by composition with $g$ preserves commutative monoid (algebra)
objects.
Thus it gives rise to $g^S:\CAlg(\uCat)\to \CAlg(\uCat)$.
Roughly speaking, $g^S$ sends a symmetric monoidal $\infty$-category
$\mathcal{C}^\otimes$
to $\Fun(S,\mathcal{C})$ endowed with the symmetric monoidal
structure
$\Fun(S,\mathcal{C})\times \Fun(S,\mathcal{C})\to \Fun(S,\mathcal{C})$
given by symmetric monoidal structure
$\mathcal{C}\times \mathcal{C}\to \mathcal{C}$.
We informally regard an object in $\Fun(S,\mathcal{C})$
as something like a fiber bundle of objects in $\mathcal{C}$ over the geometric realization $|S|$.
Let us consider the composite
\[
\Theta^S:\CAlg_R\stackrel{\Theta}{\longrightarrow} \CAlg(\uCat)\stackrel{g^S}{\longrightarrow} \CAlg(\uCat).
\]
This is a geometric $R$-system.
\end{Example}

\vspace{2mm}

{\it Automorphism group functor.}
Let $\mathcal{C}^\otimes$ be a symmetric monoidal small
$\infty$-category.
Let $\omega:\mathcal{C}^\otimes \to \mathcal{T}^\otimes(R)$
be a symmetric monoidal functor.
We write $\mathcal{C}$ for its underlying $\infty$-category.
Let $\mathcal{T}^\otimes$ be a geometric $R$-system.
Let $\theta_{\mathcal{C}^\otimes}:\CAlg(\uCat)\to \mathcal{S}$ be the functor corresponding to $\mathcal{C}^\otimes$
via the Yoneda embedding $\CAlg(\uCat)^{op}\subset \Fun(\CAlg(\uCat),\mathcal{S})$.
Then
the composite
\[
\xi:\CAlg_R \stackrel{\mathcal{T}^\otimes}{\longrightarrow} \CAlg(\uCat) \stackrel{\theta_{\mathcal{C}^\otimes}}{\longrightarrow} \mathcal{S}
\]
carries $A$ to the space equivalent to
$\Map^\otimes(\mathcal{C}^\otimes,\mathcal{T}^\otimes(A))$.
We can extends $\xi$ to
$\xi_\ast:\CAlg_R\to \mathcal{S}_{\ast}$
by using the symmetric monoidal functor $\omega$.
Here $\mathcal{S}_\ast$ denotes the $\infty$-category of pointed
spaces, that is, $\mathcal{S}_{\Delta^0/}$.
To explain this, let $\mathcal{M}\to \CAlg_R$ be a left fibration
corresponding to $\xi$.
An extension of $\xi$ to $\xi_\ast$ amounts to giving
a section $\CAlg_R\to \mathcal{M}$ of the
left fibration
$\mathcal{M}\to \CAlg_R$.
According to \cite[3.3.3.4]{HTT}
a section corresponds to an object
in the $\infty$-category $\mathcal{L}$
which is the limit of the diagram of spaces (or $\infty$-categories) given by $\xi$;
$A\mapsto \Map^\otimes(\mathcal{C}^\otimes, \mathcal{T}^\otimes(A))$.
Thus if $\lim\mathcal{T}^\otimes(A)$ denotes the limit of
$\mathcal{T}^\otimes:\CAlg_R\to \CAlg(\uCat)$, then
$\mathcal{L}$ is equivalent to $\Map^\otimes(\mathcal{C}^\otimes,\lim\mathcal{T}^\otimes(A))$ as $\infty$-categories (or equivalently spaces).
The natural functor $\mathcal{T}^\otimes(R)\to \lim\mathcal{T}^\otimes(A)$
induces $p:\Map^\otimes(\mathcal{C}^\otimes,\mathcal{T}^\otimes(R)) \to \Map^\otimes(\mathcal{C}^\otimes,\lim\mathcal{T}^\otimes(A))\simeq \mathcal{L}$.
The image $p(\omega)$ in $\mathcal{L}$ gives rise to
a section $\CAlg_R\to \mathcal{M}$.
Consequently, we have $\xi_\ast:\CAlg_R\to \mathcal{S}_\ast$
which extends $\xi$.
We define $\Aut(\omega)$
to be the composite
\[
\CAlg_R\stackrel{\xi_\ast}{\longrightarrow} \mathcal{S}_\ast\stackrel{\Omega_\ast}{\longrightarrow} \Grp(\mathcal{S}),
\]
where the second functor is the based loop functor,
and $\Grp(\mathcal{S})$ denotes the $\infty$-category of group objects in $\mathcal{S}$.
We refer to $\Aut(\omega)$ as the automorphism group functor of $\omega:\mathcal{C}^\otimes\to \TTT^\otimes(R)$.
For any $A\in \CAlg_R$, $\Aut(\omega)(A)$
is equivalent (as an object in $\mathcal{S}$)
to
the mapping space from the symmetric monoidal functor $\mathcal{C}^\otimes\to \TTT^\otimes(R)\to \TTT^\otimes(A)$ to itself in $\Map^\otimes(\mathcal{C}^\otimes,\TTT^\otimes(A))$.
We often abuse notation and write $\Aut(\omega)$ also for
the composition $\CAlg_R\to \Grp(\mathcal{S})\to \mathcal{S}$
with the forgetful functor.

Let $\Omega:\mathcal{C}\to \mathcal{T}(R)$ be the underlying functor
of $\omega$.
Let $\theta_{\mathcal{C}}:\uCat\to \mathcal{S}$ be the functor corresponding to $\mathcal{C}$
via the Yoneda embedding $\uCat^{op}\subset \Fun(\uCat,\mathcal{S})$.
Consider the composite
\[
\eta:\CAlg_R \stackrel{\mathcal{T}}{\longrightarrow} \uCat \stackrel{\theta_{\mathcal{C}}}{\longrightarrow} \mathcal{S}.
\]
As in the above case, we can extend $\eta$ to $\eta_\ast:\CAlg_R\to \mathcal{S}_\ast$ by
$\Omega:\mathcal{C}\to \mathcal{T}(R)$.
We define $\Aut(\Omega)$
to be the composite
\[
\CAlg_R\stackrel{\eta_\ast}{\longrightarrow}  \mathcal{S}_\ast\stackrel{\Omega_*}{\longrightarrow} \Grp(\mathcal{S}).
\]
We refer to $\Aut(\Omega)$ as the automorphism group functor of $\Omega:\mathcal{C}\to \TTT(R)$.
We often abuse notation and write $\Aut(\Omega)$ also for
the composite $\CAlg_R\to \Grp(\mathcal{S})\to \mathcal{S}$.

\section{Tannakization}

The goal of this Section is to prove Theorem~\ref{supermain}.

We first prove Lemmata concerning the structure of the $\infty$-category
$\uCat$.

\begin{Lemma}
\label{gengen}
Let $\mathcal{C}$ and $\mathcal{D}$ be $\infty$-categories.
Let $F:\mathcal{C}\to \mathcal{D}$ be a functor.
Then $F$ is a categorical equivalence
if and only if the composition induces equivalences
\[
f:\Map(\Delta^0,\mathcal{C})\to \Map(\Delta^0,\mathcal{D})\ \ \ \textup{and}\ \ \ 
g:\Map(\Delta^1,\mathcal{C})\to \Map(\Delta^1,\mathcal{D})
\]
in $\mathcal{S}$.
\end{Lemma}

\Proof
The part of ``only if'' is clear. We will prove the ``if'' part.
Let $\mathcal{C}^{\simeq}$ and $\mathcal{D}^{\simeq}$
be the largest Kan complexes in $\mathcal{C}$ and $\mathcal{D}$
respectively.
The equivalence of $f$ implies that the induced map
$F^{\simeq}:\mathcal{C}^{\simeq}\to \mathcal{D}^{\simeq}$
is a homotopy equivalence (or equivalently, categorical equivalence).
It follows that $F$ is essentially surjective.
Hence it suffices to show that $F$ is fully faithful.
Let $C$ and $C'$ be objects in $\mathcal{C}$.
There exists a natural equivalence
\[
\Map_{\mathcal{C}}(C,C')\simeq \Map(\Delta^1,\mathcal{C})\times_{\Map(\partial\Delta^1,\mathcal{C})}\{(C,C')\}
\]
in $\mathcal{S}$, where $\{(C,C')\}=\Delta^0\to \Map(\partial\Delta^1,\mathcal{C})$
corresponds to $C$ and $C'$.
The induced map $\Map_{\mathcal{C}}(C,C')\to \Map_{\mathcal{D}}(F(C),F(C'))$
can be identified with
\[
\Map(\Delta^1,\mathcal{C})\times_{\Map(\partial\Delta^1,\mathcal{C})}\{(C,C')\}\to \Map(\Delta^1,\mathcal{D})\times_{\Map(\partial\Delta^1,\mathcal{D})}\{(F(C),F(C'))\},
\]
which is an equivalence in $\mathcal{S}$ by our assumption.
\QED

We will construct the full subcategory $\langle \Delta^0,\Delta^1\rangle$
of $\textup{Cat}_\infty$ by the following inductive steps.
We first note that $\textup{Cat}_\infty$ is a presentable $\infty$-category
since it is (equivalent to) the simplicial nerve of the simplicial
category consisting of fibrant objects in the combinatorial model
category of small marked simplicial sets, defined in \cite[3.1.3.7]{HTT}.
Choose a regular cardinal $\kappa$ such that
$\textup{Cat}_\infty$ is $\kappa$-accessible (cf. \cite[5.4.2]{HTT})
and both $\Delta^0$ and $\Delta^1$ are $\kappa$-compact.
Let $[\Delta^0,\Delta^1]_0$ be the full subcategory of $\textup{Cat}_\infty$,
spanned by $\Delta^0$ and $\Delta^1$.
We define a transfinite sequence
\[
[\Delta^0,\Delta^1]_0\to [\Delta^0,\Delta^1]_1\to \cdots
\]
of full subcategories indexed by ordinals smaller than $\kappa$.
Supposing that $[\Delta^0,\Delta^1]_\alpha$ has been defined,
we define $[\Delta^0,\Delta^1]_{\alpha+1}$ to be the full subcategory of
$\uCat$
spanned by retracts of colimits of $\kappa$-small diagrams taking values in $[\Delta^0,\Delta^1]_\alpha$. Here colimits are taken in $\textup{Cat}_\infty$. If $\lambda$ is a limit ordinal,
$[\Delta^0,\Delta^1]_\lambda$ is defined to be $\bigcup_{\alpha<\lambda}[\Delta^0,\Delta^1]_\alpha$. We set $\langle \Delta^0,\Delta^1\rangle=\bigcup_{\alpha<\kappa}[\Delta^0,\Delta^1]_\alpha$.

\begin{Lemma}
\label{autolem}
The full subcategory $\langle \Delta^0,\Delta^1\rangle$
has $\kappa$-small colimits which are compatible with
those in $\textup{Cat}_\infty$
Moreover, it is idempotent complete.
\end{Lemma}

\Proof
Let $f:I\to \langle \Delta^0,\Delta^1\rangle$ be a functor where $I$ is
a $\kappa$-small simplicial set. We will show that
the colimit of $f$ in $\textup{Cat}_\infty$ belongs to
$\langle \Delta^0,\Delta^1\rangle$. Since $I$ is $\kappa$-small,
we have an ordinal $\tau$ smaller than $\kappa$, such that
$f$ factors through $[\Delta^0,\Delta^1]_\tau\subset \langle \Delta^0,\Delta^1\rangle$.
Then by our construction, the colimit in $\textup{Cat}_\infty$ belongs to
$[\Delta^0,\Delta^1]_{\tau+1}$.
Since $\textup{Cat}_\infty$ is idempotent complete
and $\langle \Delta^0,\Delta^1\rangle$ is closed under retracts,
$\langle \Delta^0,\Delta^1\rangle$ is idempotent complete.
\QED.

\begin{Lemma}
\label{autolem2}
The full subcategory $\langle \Delta^0,\Delta^1\rangle$
is the smallest full subcategory having the properties:
\begin{itemize}
\item it includes $\Delta^0$ and $\Delta^1$,

\item it has $\kappa$-small colimits, and the inclusion
$\langle \Delta^0,\Delta^1\rangle\to \textup{Cat}_\infty$ preserves $\kappa$-small colimits,

\item it is idempotent complete.
\end{itemize}
Moreover, the full subcategory $\langle \Delta^0,\Delta^1\rangle$
is small.
\end{Lemma}

\Proof
By Lemma~\ref{autolem}, it will suffice to prove that for each $\alpha$,
$[\Delta^0, \Delta^1]_\alpha$ is contained in the smallest
full subcategory. We proceed by transfinite induction.
The case of $\alpha=0$ is obvious.
Suppose that $[\Delta^0,\Delta^1]_\beta$ is contained in the smallest
full subcategory where $\beta<\alpha$. Then by
our construction in both
successor and limit cases, $[\Delta^0,\Delta^1]_{\alpha}$ is so.
To see the second claim, note that the (small) full subcategory
consisting of $\kappa$-compact objects in $\textup{Cat}_\infty$
is idempotent complete and admits $\kappa$-small colimits
which are compatible with those in $\textup{Cat}_\infty$.
Thus the first claim implies that it contains $\langle \Delta^0,\Delta^1\rangle$. It follows that $\langle \Delta^0,\Delta^1\rangle$ is
small.
\QED

Let $\Ind_{\kappa}(\langle \Delta^0,\Delta^1\rangle)$ be the
full subcategory of $\Fun(\langle \Delta^0,\Delta^1\rangle^{op},\mathcal{S})$
spanned by colimits of $\kappa$-filtered diagrams taking values in
$\langle \Delta^0,\Delta^1\rangle\subset \Fun(\langle \Delta^0,\Delta^1\rangle^{op},\mathcal{S})$ (see \cite[5.3.5]{HTT}).
According to Lemma~\ref{autolem} and \cite[5.5.1.1]{HTT},
$\Ind_{\kappa}(\langle \Delta^0,\Delta^1\rangle)$
is a presentable $\infty$-category.

\begin{Corollary}
\label{coincide}
The full subcategory $\langle \Delta^0,\Delta^1\rangle$
coincides with the full subcategory $\mathcal{E}$
consisting of $\kappa$-compact
objects in $\Ind_{\kappa}(\langle \Delta^0,\Delta^1\rangle)$.
\end{Corollary}

\Proof
Since $\langle \Delta^0,\Delta^1\rangle$ is idempotent complete
by Lemma~\ref{autolem},
our assertion follows from
\cite[5.4.2.4]{HTT} which says that the natural
inclusion $\langle \Delta^0,\Delta^1\rangle\to \mathcal{E}$
is idempotent completion.
\QED

\begin{Proposition}
\label{repcolim}
Let $\theta:\Ind_{\kappa}(\langle \Delta^0,\Delta^1\rangle)\to \textup{Cat}_\infty$
be a left Kan extension of $\langle \Delta^0,\Delta^1\rangle\to \textup{Cat}_\infty$ that preserves $\kappa$-filtered colimits (cf. \cite[5.3.5.10]{HTT}).
Then $\theta$ is a categorical equivalence.
\end{Proposition}

\Proof
Note that according to Lemma~\ref{autolem}
$\langle \Delta^0,\Delta^1\rangle\to \textup{Cat}_\infty$ preserves
$\kappa$-small colimits, and by Corollary~\ref{coincide}
$\langle \Delta^0,\Delta^1\rangle$
coincides with the full subcategory of $\kappa$-compact objects
in $\Ind_{\kappa}(\langle \Delta^0,\Delta^1\rangle)$.
Therefore by \cite[5.5.1.9]{HTT} $\theta$ preserves small colimits.
Note that every object in $\langle \Delta^0,\Delta^1\rangle$
is $\kappa$-compact in $\textup{Cat}_\infty$ (see the proof of Lemma~\ref{autolem2}).
Therefore invoking \cite[5.3.5.11 (1)]{HTT} we deduce that
$\theta$ is fully faithful.
By adjoint functor theorem \cite[5.5.2.9 (1)]{HTT} to $\theta$,
there exists its right adjoint $\xi:\textup{Cat}_\infty\to \Ind_{\kappa}(\langle \Delta^0,\Delta^1\rangle)$.
Let $\mathcal{C}$ be a (small) $\infty$-category.
To prove our assertion,
it suffices to show that the counit map
$\theta\circ \xi(\mathcal{C})\to \mathcal{C}$
is a categorical equivalence.
Now it can be checked by Lemma~\ref{gengen}.
\QED

Now we show that the example presented in Example~\ref{example1}
is a geometric $R$-system.

\begin{Lemma}
\label{Homrep}
Let $M,N\in \PMod_R$.
The functor $\Hom(M,N):\CAlg_R\to \mathcal{S}$ is representable by a derived affine scheme over $R$.
Moreover, $\Aut(M):\CAlg_R \to \mathcal{S}$ is representable by a derived affine scheme over $R$. Namely, the example in Example~\ref{example1} is a geometric
$R$-system.
\end{Lemma}

\Proof
Note that there exist natural equivalences
\[
\Map_{\Mod_A}(M\otimes_RA,N\otimes_RA)\simeq \Map_{\Mod_R}(M,N\otimes_RA)\simeq \Map_{\CAlg_R}(\textup{Sym}^*(M\otimes_RN^{\vee}),A)
\]
in $\mathcal{S}$,
where $N^{\vee}$ is the dual object of $N$ in $\PMod_R$,
and $\Sym^*(M\otimes_RN^\vee)$ is a free commutative $R$-ring spectrum
determined by $M\otimes_RN^\vee$.
Consequently, we conclude that $\Spec \Sym^*(M\otimes_RN^\vee)$ represents
the functor $\Hom(M,N)$.

Next consider $\Aut(M):\CAlg_R \to \mathcal{S}$.
This case follows from \cite[II, 1.2.10.1]{HAG2}.
\QED

Let $\Omega:\mathcal{C}\to \mathcal{T}(R)$ be the underlying functor
of $\omega:\mathcal{C}^\otimes\to \TTT^\otimes(R)$
in Theorem~\ref{bigmain}.
For any $A\in \CAlg_R$, we let $\Omega_A$ to be the composite
$\mathcal{C}\to \TTT(R)\to \TTT(A)$ where the second functor is induced
by $R\to A$.
Consider the functor $\Aut(\Omega):\CAlg_R\to \mathcal{S}$
given by $A\mapsto \Map_{\Map(\mathcal{C},\PMod_A)}(\Omega_A,\Omega_A)$
(see the previous Section).

\begin{Lemma}
\label{easyrep}
Suppose that $\mathcal{C}$ is equivalent to either
$\Delta^0$ or $\Delta^1$. 
Then $\Aut(\Omega)$ is representable by a derived affine scheme over $R$.
\end{Lemma}

\Proof
We first treat the case of $\Delta^0$.
Let $M=\Omega(\{0\})\in \TTT(R)$, where $\{0\}$ denotes
the object in $\Delta^0$.
In this case, $\Aut(\Omega)$ is representable by a derived affine
scheme $\Aut(M)$ over $R$ since $\TTT^\otimes$ is a geometric $R$-system.

Next we consider the case of $\mathcal{C}=\Delta^1$.
Let $M:=\Omega(\{0\})\in \TTT(R)$ and $N:=\Omega(\{1\})\in \TTT(R)$, where $\{0\}$ and $\{1\}$ denote objects in $\Delta^1$.
Then $\Aut(\Omega)$ is representable by
the fiber product of derived affine schemes
\[
\Aut(M)\times_{\Hom(M,N)}\Aut(N)
\]
where we regard $\Hom(M,N)$ as a derived affine scheme by
(A2) of the definition of geometric $R$-systems.
This completes the proof.
\QED

Using Proposition~\ref{repcolim}
we first treat the case where we do not take account into
symmetric monoidal structures.

\begin{Proposition}
\label{nonmonoreprep}
Let $\mathcal{C}$ be a small $\infty$-category.
Then $\Aut(\Omega)$ is representable by a derived affine scheme over $R$.
\end{Proposition}

\Proof
Suppose first that $\mathcal{C}$ belongs to $\langle \Delta^0,\Delta^1\rangle$.
Recall $\langle \Delta^0,\Delta^1\rangle=\bigcup_{\alpha<\kappa}[\Delta^0,\Delta^1]_\alpha$. We suppose that $\mathcal{C}$ belongs to $[\Delta^0,\Delta^1]_\alpha$.
We proceed by transfinite induction on $\alpha$.

If $\alpha=0$, then our assertion follows from Lemma~\ref{easyrep}.
Suppose that $\alpha<\lambda$ our assertion holds.
If $\lambda$ is a limit ordinal, then the case of $\lambda$ follows from
the definition of $[\Delta^0,\Delta^1]_\lambda$.
When $\lambda$ is a successor ordinal and $\tau+1=\lambda$,
$\mathcal{C}\in [\Delta^0,\Delta^1]_\lambda$ is a retract of a colimit
of a $\kappa$-small diagram taking values in $[\Delta^0,\Delta^1]_\tau$.
If $\mathcal{C}$ is a colimit
of a $\kappa$-small diagram taking values in $[\Delta^0,\Delta^1]_\tau$,
then by the inductive assumption on $[\Delta^0,\Delta^1]_\tau$,
we see that $\Aut(\Omega)$ is expressed as a limit of a $\kappa$-small
diagram of derived affine schemes (since the $\infty$-category
of derived affine  schemes admits small limits).
Indeed, suppose that
$\mathcal{C}$ is equivalent to a colimt
$\textup{colim}_{\mu\in I}\mathcal{C}_{\mu}$ 
of small $\infty$-categories $\mathcal{C}_{\mu}$
indexed by a small $\infty$-category $I$,
and our claim holds for the case of $\mathcal{C}_{\mu}$,
that is, the automorphism group functor
$\Aut(\Omega_\mu)$ of $\Omega_\mu:\mathcal{C}_\mu \to \mathcal{T}(R)$
is
representable by a derived affine scheme $G_\mu$ over $R$
(here $\Omega\simeq \lim_{\mu\in I}\Omega_\mu$).
It follows that
$\Aut(\Omega)$ is representable by a limit of derived affine schemes $G_\mu$.
If $\mathcal{C}$ is a retract of such a colimit,
then the retract is expressed as a colimit of a certain idempotent diagram
indexed by the simplicial set $\operatorname{Idem}$
(see \cite[4.4.5.4 (1)]{HTT}).
Hence our assertion holds also
for the case of retracts.
Therefore if $\mathcal{C}$ belongs to $\langle \Delta^0,\Delta^1\rangle$,
our assertion holds.

In general case, by Proposition~\ref{repcolim},
$\mathcal{C}$ can be expressed as
a colimit of a small $\kappa$-filtered diagram taking values in
$\langle \Delta^0,\Delta^1\rangle$.
It follows from $\langle \Delta^0,\Delta^1\rangle$ that in the general case $\Aut(\Omega)$ can
be written as a $\kappa$-filtered limit of derived affine schemes over $R$.
\QED

\begin{Corollary}
The functor $\Theta^{S}:\CAlg_R\to \CAlg(\uCat)$ in Example~\ref{example2}
is a geometric $R$-system.
\end{Corollary}

\Proof
Replace $\mathcal{C}$ in the proof of Proposition~\ref{nonmonoreprep}
by the Kan complex $S$. Then the proof together with Lemma~\ref{Homrep}
implies (A1). The proof of (A2) is similar.
\QED

\begin{Theorem}
\label{bigmain}
Let $\mathcal{T}^\otimes:\CAlg_R\to \CAlg(\uCat)$ be a geometric $R$-system.
Let $\mathcal{C}^\otimes$ be a symmetric monoidal small
$\infty$-category and let
$\omega:\mathcal{C}^\otimes\to \mathcal{T}^\otimes(R)$
be a symmetric monoidal functor.
Then $\Aut(\omega)$ is representable by a derived affine group scheme
over $R$.
\end{Theorem}

\Proof
For ease of notation, we let $\Gamma=\NNNN(\FIN)$.
Note first that
a symmetric monoidal $\infty$-category can be regarded as
a commutative monoid object in $\uCat$ (see \cite[2.4.2]{HA}).
Let $\mathcal{C}^\otimes$ and $\TTT^\otimes(A)$ be symmetric monoidal
$\infty$-categories.
Hence we regard them as commutative
monoid objects
$p:\Gamma \to \textup{Cat}_\infty$
and $q_A:\Gamma \to \textup{Cat}_\infty$ respectively.
We remark that $p(\langle n\rangle_\ast)\simeq \mathcal{C}^{\times n}$ and $q_A(\langle n\rangle_\ast)\simeq \TTT(A)^{\times n}$.
We let $r_A:\Gamma\times \Delta^1\to \uCat$ be the map corresponding
to the composite $\omega_A:\mathcal{C}^\otimes \to \TTT^\otimes(R)\to \TTT^\otimes(A)$.
Then by using \cite[4.2.1.8]{HTT} twice $\Aut(\omega)(A)$ can be identified with the Kan complex
\[
\Fun(\Gamma\times\Delta^1\times\Delta^1,\uCat)\times_{\Fun(\Gamma\times \partial(\Delta^1\times \Delta^1),\uCat)}\{(c_p\sqcup c_{q_A})\cup(r_A\sqcup r_A)\},
\]
where $\{(c_p\sqcup c_{q_A})\cup(r_A\sqcup r_A)\}$ denotes the union
$(c_p\sqcup c_{q_A})\cup(r_A\sqcup r_A):\Gamma\times \partial(\Delta^1\times\Delta^1) \to \uCat$
such that $c_p:\Gamma \times\Delta^1\times\{0\} \stackrel{\textup{pr}_1}{\to} \Gamma \stackrel{p}{\to} \uCat$,
$c_{q_A}:\Gamma \times\Delta^1\times\{1\} \stackrel{\textup{pr}_1}{\to} \Gamma \stackrel{q_A}{\to} \uCat$, and
$r_A\sqcup r_A:\Gamma \times \partial\Delta^1\times \Delta^1\to \uCat$.
Thus $\Aut(\omega)$
is given by
\[
A \mapsto \Map(\Gamma\times \Delta^1\times\Delta^1, \uCat)\times_{\Map(\Gamma\times \partial(\Delta^1\times \Delta^1), \uCat)}\{(c_p\sqcup c_{q_A})\cup(r_A\sqcup r_A)\},
\]
where the right hand side is the (homotopy) fiber product.
Replacing $\kappa$ above by a larger regular cardinal if
necessary (cf. Proposition~\ref{repcolim}),
we may assume that $\Gamma$ is $\kappa$-compact in $\uCat$.
Let $f:I\to \Gamma$ be a functor
from $I\in [\Delta^0,\Delta^1]_\alpha$ to
$\Gamma$.
Consider the composite
\[
\CAlg_R\stackrel{\TTT^\otimes}{\longrightarrow} \Fun(\Gamma,\uCat)\to \Fun(I,\uCat)\to \mathcal{S}
\] 
where the first functor is $\TTT^\otimes:\CAlg_R\to \CAlg(\uCat)\subset \Fun(\Gamma,\uCat)$, and the second functor is induced by the composition with
$f$, and the third functor is representable by $p \circ f$.
By $f^*r_R:=r_R\circ (f\times \textup{Id}_{\Delta^1}):I\times \Delta^1\to \uCat$, 
we can extends the above composite to $\CAlg_R\to \mathcal{S}_\ast$
as in the previous Section.
Composing with $\mathcal{S}_\ast\stackrel{\Omega_*}{\to} \Grp(\mathcal{S})\to \mathcal{S}$
we have $\CAlg_R\to \mathcal{S}$, which we shall denote by $\Aut(\omega)_f$.
This functor sends $A$ to the (homotopy) fiber product
\[
\Map(I\times \Delta^1\times \Delta^1, \uCat)\times_{\Map(I\times \partial(\Delta^1\times \Delta^1), \uCat)}\{(f^*c_p\sqcup f^*c_{q_A})\cup(f^*r_A\sqcup f^*r_A)\}.
\]
We claim that if $I$ is $\kappa$-compact then $\Aut(\omega)_f$
is representable by a derived affine scheme over $R$.
Suppose that $I$ belongs to $[\Delta^0,\Delta^1]_0$.
Then the case of $I\simeq \Delta^0$ is reduced to
Proposition~\ref{nonmonoreprep};
suppose that the image $f(\Delta^0)$ corresponds to $\langle n\rangle_\ast$.
Recall that $q_A(\langle n\rangle_*)$ is equivalent to the $n$-fold product $\TTT(A)^{\times n}$ as $\infty$-categories.
In this case, $\Aut(\omega)_f$ is
given by $\Aut(\omega_n):\CAlg_R\to \mathcal{S}$,
\[
A\mapsto \prod_{1\le i\le n}\Aut(\textup{pr}_i\circ \omega_{n,A})
\]
where $\omega_{n,A}$ is the functor $p(\langle n\rangle_*)\to q_A(\langle n\rangle_*)$ induced by $\omega$, and $\textup{pr}_i:\TTT(A)^{\times n}\to \TTT(A)$
is the $i$-th projection.
 Hence thanks to Proposition~\ref{nonmonoreprep},
this functor $\Aut(\textup{pr}_i\circ \omega_{n,A})$
is representable by a derived affine scheme over $R$.
It follows that $\Aut(\omega_n)$ is representable by a derived affine scheme
over $R$.
When $f:I\simeq \Delta^1$ and $I\to \Gamma$
corresponds to
$\langle m \rangle_*\to \langle n\rangle_*$,
$\Aut(\omega)_f$
is representable by
\[
\Aut(\omega_m)\times_{\Aut(\omega_{m,n})}\Aut(\omega_n)
\]
where $\omega_{m,n}$ is the functor $p(\langle m\rangle_*)\to q(\langle n\rangle_*)$ induced by
$\omega$ and $\langle m\rangle_*\to \langle n\rangle_*$. Thus this case is again reduced to Proposition~\ref{nonmonoreprep}.
Next suppose that if $\alpha<\lambda$
our assertion holds for $\alpha$.
If $\lambda$ is a limit ordinal, our assertion also holds for the
case of $\lambda$.
Assume that $\lambda$ is a successor ordinal and $\tau+1=\lambda$.
Let $I\to \Gamma$ be a functor with $I\in [\Delta^0,\Delta^1]_\lambda$
and consider the case when $I\simeq \textup{colim}I_\mu$,
where $\textup{colim}I_\mu$ is a colimit of a
$\kappa$-small diagram taking values in $[\Delta^0,\Delta^1]_\tau$.
According to \cite[1.2.13.8]{HTT}, $I\simeq \textup{colim}I_\mu\to \Gamma$
is a colimit also in $(\uCat)_{/\Gamma}$.
Note that the cartesian product commutes with
colimits in $\uCat$.
Thus the assumption for the case of $\tau$ (and the definition
of $\Aut(\omega)_f$) implies that
our assertion also holds for the
case of $\lambda$.
If $I'$ is a retract of the above $I$,
a retract can also be expressed as the colimit (see \cite[4.4.5]{HTT}).
Hence our assertion holds for the case of the retract.
This implies that for every $\kappa$-compact $\infty$-category $I$,
our assertion holds. In particular, if $I=\Gamma$,
Theorem~\ref{bigmain} follows since $\Aut(\omega):\CAlg_R\to \Grp(\mathcal{S})\to \mathcal{S}$ is representable by a derived affine scheme over $R$.
\QED

\begin{Proposition}
\label{mainreprep}
Let $\omega:\mathcal{C}^\otimes\to \PMod_R^\otimes$ be
a symmetric monoidal functor where $\mathcal{C}^\otimes$
is a symmetric monoidal small $\infty$-category.
(Here the geometric $R$-system is given in Example~\ref{example1}.)
Then the functor $\Aut(\omega)$ is representable by a derived affine
group scheme $G$ over $R$.
\end{Proposition}

\Proof
It follows from Theorem~\ref{bigmain} and Lemma~\ref{Homrep}.
\QED

Let $\mathcal{C}^\otimes$ be a symmetric monoidal small $\infty$-category
and let $\omega:\mathcal{C}^\otimes\to \PMod_R^\otimes$ be a
symmetric monoidal functor.
Let $G$ be a derived affine group scheme
over $R$.
Let $\PRep_G^\otimes$ be the symmetric monoidal
stable $\infty$-category of perfect representations of $G$
(see A.6).
Suppose that $\omega$ is extended to
a symmetric monoidal functor
$\mathcal{C}^\otimes\to \PRep_{G}^\otimes$.
Namely, the composite $\mathcal{C}^\otimes\to \PRep_G^\otimes\to \PMod_R^\otimes$ with the forgetful functor is equivalent to $\omega$.
Next our goal is Proposition~\ref{autotan} which relates such extensions
with actions on $\omega$.
Let $\NNNN(\Delta)^{op}\to \Aff_R\subset \Fun(\CAlg_R,\widehat{\mathcal{S}})$
be a functor corresponding to $G$ and let
$\mathsf{B}G$ be its colimit.
Let $(\Aff_R)_{/\mathsf{B}G}$ be the full subcategory of $\Fun(\CAlg_R,\widehat{\mathcal{S}})_{/\mathsf{B}G}$ spanned by objects $X\to \mathsf{B}G$
such that $X$ are affine schemes, that is, objects which belong to the essential image of Yoneda embedding $\Aff_R\hookrightarrow \Fun(\CAlg_R,\widehat{\mathcal{S}})$. There is the natural projection $(\Aff_R)_{/\mathsf{B}G}\to \Aff_R$, that is a right fibration (cf. \cite[2.0.0.3]{HTT}).
Let $\pi:\Spec R\to \mathsf{B}G$ be the natural projection.
This determines a map between right fibrations
$\Aff_R=(\Aff_R)_{/\Spec R}\to (\Aff_R)_{/\mathsf{B}G}$
over $\Aff_R$.
Let $(\Aff_R)_{/\mathsf{B}G} \to \mathcal{S}^{op}$
be a functor which assigns $\Map^\otimes(\mathcal{C}^\otimes,\PMod^\otimes_A)$ to $\Spec A$ in $(\Aff_R)_{/\mathsf{B}G}$.
Here $\Map^\otimes(-,-)$ indicates the mapping space
in $\CAlg(\uCat)$.
More precisely, let 
\[
c:(\Aff_R)_{/\mathsf{B}G} \to \Aff_R \stackrel{\theta}{\to} \CAlg(\uCat)^{op} \to \mathcal{S}^{op}
\]
be the composition where the first functor is the natural projection,
and the third is the image of $\mathcal{C}^\otimes$
by Yoneda embedding $(\CAlg(\uCat))^{op}\to \Fun(\CAlg(\uCat),\mathcal{S})$.
Let $\theta:\Aff_R\to \CAlg(\uCat)^{op}$
be the functor induced by $\Theta$, which carries $\Spec A$ to $\PMod_A^\otimes$.
By the unstraightening functor \cite[3.2]{HTT} together
with \cite[4.2.4.4]{HTT}
the composition $(\Aff_R)_{/\mathsf{B}G}\to \mathcal{S}^{op}$
gives rise to a right fibration $p:\mathcal{M}\to (\Aff_R)_{/\mathsf{B}G}$.
The mapping space
$\Map^\otimes(\mathcal{C}^\otimes,\PRep_G^\otimes)$
is homotopy equivalent to the limit of spaces
\[
\lim_{\Spec A\to \mathsf{B}G}\Map^\otimes(\mathcal{C}^\otimes, \theta(\Spec A))
\]
where $\Spec A\to \mathsf{B}G$ run over $(\Aff_R)_{/\mathsf{B}G}$
and $\PMod_{\mathsf{B}G}^\otimes\simeq \lim_{\Spec A\to \mathsf{B}G}\theta(\Spec A)$ (see A.6 for $\PMod_{\mathsf{B}G}^\otimes$).

\begin{Lemma}
If we denote by
$\Map_{(\Aff_R)_{/\mathsf{B}G}}((\Aff_R)_{/\mathsf{B}G},\mathcal{M})$
the simplicial set of the sections of
$p:\mathcal{M}\to (\Aff_R)_{/\mathsf{B}G}$
(i.e., the set of $n$-simplexes of $\Map_{(\Aff_R)_{/\mathsf{B}G}}((\Aff_R)_{/\mathsf{B}G},\mathcal{M})$ is the set of $(\Aff_R)_{/\mathsf{B}G}\times \Delta^n\to \mathcal{M}$ over $(\Aff_R)_{/\mathsf{B}G}$),
then
there is a categorical equivalence
$\Map^\otimes(\mathcal{C}^\otimes,\PMod_{\mathsf{B}G}^\otimes) \simeq \Fun_{(\Aff_R)_{/\mathsf{B}G}}((\Aff_R)_{/\mathsf{B}G},\mathcal{M})$.
\end{Lemma}

\Proof
It follows from \cite[3.3.3.2]{HTT}.
\QED

\begin{Proposition}
\label{autotan}
There is a natural equivalence
\[
\Map_{\CAlg(\uCat)_{/\PMod_R^\otimes}}(\mathcal{C}^\otimes,\PRep_{G}^\otimes) \simeq \Map_{\Fun(\CAlg_R,\Grp(\mathcal{S}))}(G,\Aut(\omega))
\]
in $\mathcal{S}$.
This equivalence is functorial in the following sense:
Let $L:\textup{GAff}_R \to \mathcal{S}^{op}$ be the functor
which assigns $G$ to $\Map_{\CAlg(\uCat)_{/\PMod_R^\otimes}}(\mathcal{C}^\otimes,\PRep_{G}^\otimes)$.
Let $M:\textup{GAff}_R \to \mathcal{S}^{op}$ be the functor which assigns
$G$ to
$\Map_{\Fun(\CAlg_R,\Grp(\mathcal{S}))}(G,\Aut(F))$.
(See the proof below for the formulations of $L$ and $M$.)
Then there exists a natural equivalence from $L$ to $M$.
\end{Proposition}

\Proof
In order to make our proof readable we first show the first assertion without defining $L$ and $M$.
The mapping space
$\Map_{\CAlg(\uCat)_{/\PMod_R^\otimes}}(\mathcal{C}^\otimes,\PMod_{\mathsf{B}G}^\otimes)$
is the homotopy limit (i.e. the limit in $\mathcal{S}$)
\[
\Map^\otimes(\mathcal{C}^\otimes,\PMod_{\mathsf{B}G}^\otimes)\times_{\Map^\otimes(\mathcal{C}^\otimes,\PMod_R^\otimes)}\{\omega\}
\]
where $\{\omega\}=\Delta^0\to \Map^\otimes(\mathcal{C}^\otimes,\PMod_R^\otimes)$
is determined by $\omega$.
The fiber product of Kan complexes
\[
P=\Map_{(\Aff_R)_{/\mathsf{B}G}}((\Aff_R)_{/\mathsf{B}G},\mathcal{M}) \times_{\Map_{(\Aff_R)_{/\mathsf{B}G}}(\Aff_R,\mathcal{M})}\{\omega\}
\]
is a homotopy limit since $\Aff_R\to (\Aff_R)_{/\mathsf{B}G}$
is a monomorphism (that is, a cofibration in the Cartesian simplicial
model category of (not necessarily small) marked simplicial sets
$(\widehat{\textup{Set}}^+_{\Delta})_{/(\Aff_R)_{/\mathsf{B}G}}$, see \cite[3.1.3.7]{HTT})
and thus the induced map is a Kan fibration.
Here $\Delta^0=\{\omega\}\to \Map_{(\Aff_R)_{/\mathsf{B}G}}(\Aff_R,\mathcal{M})$ is determined by $\omega$.
Let $\mathcal{N}:=\mathcal{M}\times_{(\Aff_R)_{/\mathsf{B}G}}\Aff_R$
where $\Aff_R\to (\Aff_R)_{/\mathsf{B}G}$ is determined by
the natural map
$\Spec R\to \mathsf{B}G$.
Using the Cartesian equivalence $\mathcal{N}\times_{\Aff_R}(\Aff_R)_{/\mathsf{B}G} \simeq \mathcal{M}$ over $(\Aff_{R})_{/\mathsf{B}G}$
we have homotopy equivalences
\[
\Map_{(\Aff_R)_{/\mathsf{B}G}}((\Aff_R)_{/\mathsf{B}G},\mathcal{M}) \simeq \Map_{\Aff_R}((\Aff_R)_{/\mathsf{B}G},\mathcal{N})
\]
and
\[
\Map_{(\Aff_R)_{/\mathsf{B}G}}(\Aff_R,\mathcal{M})\simeq \Map_{\Aff_R}(\Aff_R,\mathcal{N}).
\]
Thus 
$P$ is homotopy equivalent to the fiber product
\[
Q=\Map_{\Aff_R}((\Aff_R)_{/\mathsf{B}G},\mathcal{N})\times_{\Map_{\Aff_R}(\Aff_R,\mathcal{N})}\{\omega\}
\]
which is also a homotopy limit, where $\Delta^0=\{\omega\}\to  \Map_{\Aff_R}(\Aff_R,\mathcal{N})$ is determined by the section
$\Aff_R\to \mathcal{N}$ corresponding to $\omega:\mathcal{C}^\otimes \to \PMod_R^\otimes$.
We let $\alpha_{\mathsf{B}G}:\CAlg_R\to \mathcal{S}$ be a map
correspoindig to the right fibration
$(\Aff_R)_{/\mathsf{B}G}\to \Aff_R$
via the straightening functor.
There is the natural transformation
$\alpha_\ast\to \alpha_{\mathsf{B}G}$ determined by
$\Aff_R\to (\Aff_R)_{/\mathsf{B}G}$, which
we consider to be a functor $\CAlg_R\to \mathcal{S}_{\ast,\ge1}$.
Here $\mathcal{S}_{\ast,\ge1}$ denotes the full subcategory of
$\mathcal{S}_\ast$ spanned by pointed connected spaces.
Let $\alpha_{\mathcal{N}}:\CAlg_R\to \mathcal{S}_\ast$
be a functor corresponding to the right fibration
$\mathcal{N}\to \Aff_R$ equipped with
the section $\Aff_R\to \mathcal{N}$.
Observe that $\Map_{\Fun(\CAlg,\mathcal{S}_\ast)}(\alpha_{\mathsf{B}G},\alpha_{\mathcal{N}})$ is homotopy equivalent to $Q$.
By composition with $\Omega_\ast:\mathcal{S}_{\ast}\to \Grp(\mathcal{S})$
we have $G:\CAlg_R \stackrel{\mathsf{B}G}{\to} \mathcal{S}_{\ast,\ge1}\simeq \Grp(\mathcal{S})$ (that is, the composition is the original derived group scheme
$G$). Let $\alpha'_{\mathcal{N}}$ be an object in $\Fun(\CAlg_R,\mathcal{S}_{\ast,\ge 1})$ such that
$\alpha'_{\mathcal{N}}(A)$ is the pointed connected component
determined by $\alpha_{\mathcal{N}}(A)$.
Then we obtain 
\begin{eqnarray*}
Q &\simeq& \Map_{\Fun(\CAlg_R,\mathcal{S}_\ast)}(\alpha_{\mathsf{B}G},\alpha_{\mathcal{N}}) \\
&\simeq& \Map_{\Fun(\CAlg_R,\mathcal{S}_{\ast})}(\alpha_{\mathsf{B}G},\alpha'_{\mathcal{N}}) \\
&\simeq& 
\Map_{\Fun(\CAlg_R,\Grp(\mathcal{S}))}(G,\Aut(\omega)).
\end{eqnarray*}

Next to see (and formulate) the latter assertion, we will define $L$ and $M$.
We first define $L$. Since a derived affine group scheme is a group object
in the Cartesian symmetric monoidal $\infty$-category of $\Aff_R$,
thus $\textup{GAff}_R$ is naturally embedded into
$\Fun(\NNNN(\Delta)^{op},\Fun(\CAlg_R,\mathcal{S}))$
as a full subcategory.
Let $\Fun(\NNNN(\Delta)^{op},\Fun(\CAlg_R,\mathcal{S}))
\to \Fun(\CAlg_R,\mathcal{S})$
be the functor taking each simplicial object
$\NNNN(\Delta)^{op}\to \Fun(\CAlg_R,\widehat{\mathcal{S}})$
to its colimit.
Let $\rho:\textup{GAff}_R \to \Fun(\CAlg_R,\mathcal{S})$
be the composition.
Note that $G$ maps to $\mathsf{B}G$.
By the straightening and unstraightening functors \cite[3.2]{HTT}
together with
\cite[4.2.4.4]{HTT},
we have the categorical equivalence
$\Fun(\CAlg_R,\widehat{\textup{Cat}}_\infty)\simeq \NNNN(((\widehat{\textup{Set}}_\Delta^{+})_{/\Aff_R})^{cf})$ where $(\widehat{\textup{Set}}_\Delta^{+})_{/\Aff_R}$
is the category of (not necessarily small) marked simplicial sets,
which is endowed with
the Cartesian model structure in \cite[3.1.3.7]{HTT}
and $(-)^{cf}$ indicates full simplicial subcategory of cofibrant-fibrant
objects.
In particular, there is the fully faithful
functor
$\Fun(\CAlg_R,\widehat{\mathcal{S}})\to \NNNN(((\widehat{\textup{Set}}_\Delta^{+})_{/\Aff_R})^{cf})$ which carries $\mathsf{B}G$ to $(\Aff_R)_{/\mathsf{B}G}\to \Aff_R$.
Composing all these functors
 we have the composition
\[
\textup{GAff}_R \stackrel{\rho}{\to} \Fun(\CAlg_R,\mathcal{S}) \to \NNNN(((\widehat{\textup{Set}}_\Delta^{+})_{/\Aff_R})^{cf}).
\]
Since $\textup{GAff}_R\simeq (\textup{GAff}_R)_{\Spec R/}$,
the composition is extended to
$u:\textup{GAff}_R\to \NNNN(((\widehat{\textup{Set}}_\Delta^{+})_{/\Aff_R})^{cf})_{\Aff_R/}$.
Through Yoneda embedding \[\NNNN(((\widehat{\textup{Set}}_\Delta^{+})_{/\Aff_R})^{cf})_{\Aff_R/}\to \Fun((\NNNN(((\widehat{\textup{Set}}_\Delta^{+})_{/\Aff_R})^{cf})_{\Aff_R/})^{op},\widehat{\mathcal{S}})\]
we define $I:(\NNNN(((\widehat{\textup{Set}}_\Delta^{+})_{/\Aff_R})^{cf})_{\Aff_R/})^{op}\to\widehat{\mathcal{S}}$ to be
the functor corresponding to $\mathcal{N}\to \Aff_R$ equipped with the section $\omega$.
Composing $I^{op}$ with
$\textup{GAff}_R\to \NNNN(((\widehat{\textup{Set}}_\Delta^{+})_{/\Aff_R})^{cf})_{\Aff_R/}$
we define $L$ to be $\textup{GAff}_R\to \widehat{\mathcal{S}}^{op}$.
To define $M$, consider the functor $\Fun(\CAlg_R,\Grp(\mathcal{S}))\to \widehat{\mathcal{S}}^{op}$ determined by $\Aut(\omega)$ via Yoneda embedding.
Then we define $M$ to be the composition
\[
\textup{GAff}_R\hookrightarrow \Fun(\CAlg_R,\Grp(\mathcal{S}))\to \widehat{\mathcal{S}}^{op}.
\]

To obtain $L\simeq M$,
note that the unstraightening functor induces a fully faithful functor
$\Fun(\CAlg_R,\widehat{\mathcal{S}}_\ast) \subset \NNNN(((\widehat{\textup{Set}}_\Delta^{+})_{/\Aff_R})^{cf})_{\Aff_R/}$.
Let $\mathsf{N}:\CAlg_R\to \mathcal{S}_\ast$ be a functor
corresponding to $\mathcal{N}\to \Aff_R$ equipped with the section $\omega$,
that is, $\mathsf{N}$ corresponds to $\alpha_\ast \to \alpha_{\mathcal{N}}$.
Let $\Fun(\CAlg_R, \mathcal{S}_\ast)\to \widehat{\mathcal{S}}^{op}$
be the functor determined by $\mathsf{N}$ via Yoneda embedding.
The functor $L$ is equivalent to
\[
\textup{GAff}_R \stackrel{u}{\to} \Fun(\CAlg_R,\widehat{\mathcal{S}}_\ast)\subset \NNNN(((\widehat{\textup{Set}}_\Delta^{+})_{/\Aff_R})^{cf})_{\Aff_R/}\to \widehat{\mathcal{S}}^{op}.
\]
Since the essential image of $\textup{GAff}_R$ in $\Fun(\CAlg_R,\mathcal{S}_{\ast})$
is contained in $\Fun(\CAlg_R,\mathcal{S}_{\ast,\ge1})$, for our purpose
we may and will replace $\alpha_{\mathcal{N}}$
by $\alpha_{\mathcal{N}}'$ (in the construction of $\mathsf{N}$)
and assume that $\mathsf{N}$ belongs to $\Fun(\CAlg_R,\mathcal{S}_{\ast,\ge1})$.
Then we see that $L$ is equivalent to
\[
\textup{GAff}_R\to  \Fun(\CAlg_R,\mathcal{S}_{\ast,\ge1})\simeq \Fun(\CAlg_R,\Grp(\mathcal{S}))\to \widehat{\mathcal{S}}^{op}
\]
where the first functor is induced by $u$ and the third functor
is determined by $\Aut(\omega)$ via Yoneda embedding.
Now the last composition is equivalent to $M$.
\QED

Now we are ready to prove the following:

\begin{Theorem}
\label{supermain}
There are a derived affine
group scheme $G$ over $R$
and a symmetric monoidal functor
$u:\mathcal{C}^\otimes\to \PRep_G^\otimes$ which makes the outer triangle in
\[
\xymatrix{
     &  \PMod_G^\otimes \ar@{..>}[d] \ar[ddr]^{\textup{forget}}  &     \\
   & \PRep_H^\otimes \ar[rd]_{\textup{forget}} &  \\
\mathcal{C}^\otimes \ar[ru] \ar[ruu]^{u} \ar[rr]_\omega & & \PMod_R^\otimes
}
\]
commute in $\CAlg(\uCat)$
such that these possess the following universality:
for any inner triangle consisting of solid arrows in the above diagram
where $H$ is a derived affine group scheme over $R$,
there exists a
morphism $f:H\to G$ of derived affine group
schemes which induces $\PRep_G^\otimes\to \PRep_H^\otimes$
(indicated by the dotted arrow) filling the above diagram.
Such $f$ is unique up to a contractible space of choices.
Moreover, the automorphism group functor $\Aut(\omega)$ is represented by $G$.
\end{Theorem}

\Proof
Take a derived affine group scheme $G$ over $R$
which represents $\Aut(\omega)$ by Proposition~\ref{mainreprep}.
By Proposition~\ref{autotan}, 
we have a symmetric monoidal functor
$\mathcal{C}^\otimes\to \PRep^\otimes_G$
that corresponds to the identity $G\simeq \Aut(\omega)\to \Aut(\omega)$.
Then Proposition~\ref{autotan} implies our claim.
\QED

We usually refer to $(G,\mathcal{C}^\otimes\stackrel{u}{\to} \PMod_G^\otimes)$
(or simply $G$)
in Theorem~\ref{supermain}
as the tannakization of $\omega:\mathcal{C}^\otimes\to \PMod_R^\otimes$.

\vspace{1mm}

The following properties are easy but useful.

\begin{Proposition}
\label{colicoli}
Let $\{\mathcal{C}^\otimes_i\}_{i\in I}$ be a
(small) collection of symmetric monoidal ful subcategories of $\mathcal{C}^\otimes$. Assume that for any finite subset $J\subset I$,
there is some $i\in I$ such that
$\bigcup_{j\in J}\mathcal{C}_i\subset \mathcal{C}_i$.
Suppose further that $\bigcup_{i\in I}\mathcal{C}_i=\mathcal{C}$.
Let $\omega_i:\mathcal{C}^\otimes_i\hookrightarrow \mathcal{C}^\otimes\stackrel{\omega}{\to} \PMod_R^\otimes$ be the composite and let $G_i$ be
the tannakization of the composite.
Then if $G$ denotes the tannakization of $\omega$, then
$G\simeq \lim_{i\in I} G_i$.
\end{Proposition}

\Proof
The collection
$\{\mathcal{C}^\otimes_i\}_{i\in I}$
constitutes a filtered partially ordered set ordered by inclusions.
As a consequence, according to \cite[3.2.3.2]{HA},
the condition
$\bigcup_{i\in I}\mathcal{C}_i=\mathcal{C}$
implies that $\mathcal{C}^\otimes$ is
a colimit of $\{\mathcal{C}^\otimes_i\}_{i\in I}$
in $\CAlg(\uCat)$.
It implies our claim (by noting the limit of derived
affine schemes commutes with the limit as functors $\CAlg_R\to \mathcal{S}$).
\QED

\begin{Proposition}
We adopt the notion of the previous Proposition.
Let $R\to R'$ be a morphism in $\CAlg$.
Then the tannakization of the composite $\mathcal{C}^\otimes\stackrel{\omega}{\to} \PMod_R^\otimes\stackrel{\otimes_RR'}{\to} \PMod_{R'}^\otimes$
is $G\times_{\Spec R}\Spec R'$.

\end{Proposition}

\section{Derived motivic Galois group}

In this Section we will
construct derived motivic Galois groups of mixed motives,
their variants, and truncated (underived) motivic Galois groups.
The term ``derived'' in the title of this Section
stems from the tannakization of the ``highly structured'' category:
stable $\infty$-category of mixed motives (see Remark~\ref{abelianmotives}).
For our purposes, we
apply Theorem~\ref{supermain} to
the stable $\infty$-category of mixed motives endowed with
the homological realization functor of a mixed Weil cohomology.
To this end, we need to construct
the realization functor of a mixed Weil cohomology theory
in the $\infty$-categorical setting.

\subsection{$\infty$-category of mixed motives}
We construct the $\infty$-category of mixed motives.
We first construct a stable $\infty$-category of motivic spectra.
There are several approaches to construct it.
Let $S$ be a scheme separated and of finite type over $\mathbb{Z}$.
Let $\SM$ be the category of smooth scheme
separated and of finite type over $S$.
One can perform the construction of Morel and Voevodsky
(\cite{MV}, \cite{V}) in the setting of $\infty$-categories.
On the other hand,
there are several model categories of motivic spectra
(e.g., \cite{Jar}, \cite{Ho2}, \cite{DR}, \cite{CD2}).
Then the passage from model categories to $\infty$-categories
allows us to have an $\infty$-category of motivic spectra.
In this paper we will adopt the latter approach.
Especially, we use the model category of symmetric Tate
spectra described in
\cite[1.4.3]{CD2},
where Cisinski and D\'eglise introduced
the theory of the mixed Weil theory which gives us the very
powerful method for constructing realization functors.

\vspace{2mm}

{\it Symmetric Tate spectra.}
We shall refer ourselves to \cite{CD1} and \cite{CD2}
for the model category of symmetric Tate spectra.
We here recall the minimal definitions for symmetric 
Tate spectra.
Let $R$ be an (ordinary) commutative ring
and $\operatorname{Sh}(\SM,R)$ the abelian category
of Nisnevich sheaves of $R$-modules.
Let $\Comp(\operatorname{Sh}(\SM,R))$
be the category of complexes of objects in $\operatorname{Sh}(\SM,R)$.
This is a symmetric monoidal category.
For the symmetric monoidal structure of complexes
of objects in a symmetric monoidal abelian category, see e.g. \cite[3.1]{CD1}.
For any $X\in \SM$, we write $R(X)$ for the Nisnevich sheaf associated to
the presheaf given by
$Y\mapsto \oplus_{f\in \Hom_{\SM}(Y,X)}R\cdot f$
where $\oplus_{f\in \Hom_{\SM}(Y,X)}R\cdot f$
is the free $R$-module generated by the set $\Hom_{\SM}(Y,X)$.
It gives rise to a functor $\SM\to \Comp(\operatorname{Sh}(\SM,R))$.
Let $R(1)[1]\in \Comp(\operatorname{Sh}(\SM,R))$
be the cokernel of the split monomorphism
$R(S)\to R(\mathbb{G}_m)$ determined by the unit $S\to \mathbb{G}_m=\Spec S[t,t^{-1}]$
of the torus.
A symmetric Tate sequence is
a sequence $\{E_n\}_{n\in \mathbb{N}}$
where $E_n$ is an object of $\Comp(\operatorname{Sh}(\SM,R))$
which is equipped with an action by the symmteric group
$\mathfrak{S}_n$
for each $n\in \mathbb{N}$.
A morphism $\{E_n\}_{n\in \mathbb{N}}\to \{F_n\}_{n\in \mathbb{N}}$
is a collection of $\mathfrak{S}_n$-equivariant maps $E_n\to F_n$.
Let $\operatorname{S}_{\textup{Tate}}(R)$ be the category of
symmetric Tate sequences.
Let $\mathfrak{S}'$ be the category of finite sets whose morphisms
are bijections.
Then the category of functors from $\mathfrak{S}'$ to $\Comp(\operatorname{Sh}(\SM,R))$ is naturally equivalent to the category of symmetric Tate
sequences (To $F:\mathfrak{S}'\to \Comp(\operatorname{Sh}(\SM,R))$
we associate $\{E_n=F(\bar{n})\}_{n\in \mathbb{N}}$ if $\bar{n}$ is $\{1,\ldots,n\}$).
For $E,F:\mathfrak{S}'\to \Comp(\operatorname{Sh}(\SM,R))$,
the tensor product
is defined to be
$\mathfrak{S}'\to \Comp(\operatorname{Sh}(\SM,R))$
given by $N\mapsto \bigoplus_{N=P\sqcup Q}E(P)\otimes F(Q)$.
It yields a symmetric monoidal structure on the category of
symmetric Tate sequences.
Let $\operatorname{Sym}(R(1))$ denote a symmetric Tate sequence
$\{R(1)^{\otimes n}\}_{n\in \mathbb{N}}$
such that $\mathfrak{S}_n$ acts on $R(1)^{\otimes n }$ by permutation.
The canonical
isomorphism $R(1)^{\otimes n}\otimes R(1)^{\otimes m}\to R(1)^{\otimes n+m}$
is $\mathfrak{S}_n\times \mathfrak{S}_m$-equivariant
when $\mathfrak{S}_n\times \mathfrak{S}_m$ acts on $R(1)^{\otimes n+m}$
through the natural inclusion $\mathfrak{S}_n\times \mathfrak{S}_m\to \mathfrak{S}_{n+m}$.
Unwinding the definition of tensor product of symmetric Tate sequences
we have
a morphism
\[
\operatorname{Sym}(R(1))\otimes \operatorname{Sym}(R(1))\to \operatorname{Sym}(R(1))
\]
which makes $\operatorname{Sym}(R(1))$ a commutative algebra object
in $\operatorname{S}_{\textup{Tate}}$.
Let $\operatorname{Sp}_{\textup{Tate}}(R)$ be the category
of modules in $\operatorname{S}_{\textup{Tate}}(R)$
over the commutative algebra object
$\operatorname{Sym}(R(1))$.
We call an object in
$\operatorname{Sp}_{\textup{Tate}}(R)$
a symmetric Tate spectrum.
In \cite[1.4.2]{CD2}, the classes of
stable $\mathbb{A}^1$-equivalences, stable $\mathbb{A}^1$-fibrations
are defined (these are important,
but we will not recall them here since we need preliminaries).
In \cite[1.4.3]{CD2} (see also \cite{CD1}), the model category structure 
of $\operatorname{Sp}_{\textup{Tate}}(R)$ is constructed:

\begin{Proposition}
The category $\operatorname{Sp}_{\textup{Tate}}(R)$
is a stable proper cellular symmetric monoidal model
category
with stable $\mathbb{A}^1$-equivalences as weak equivalences,
and stable $\mathbb{A}^1$-fibrations as fibrations.
\end{Proposition}

\begin{Remark}
A pointed model category is stable if the suspention functor
induces an equivalence of the homotopy category
(cf. \cite{Ho1}).
\end{Remark}

\begin{Lemma}
The category $\operatorname{Sp}_{\textup{Tate}}(R)$
is presentable.
In particular, it is a combinatorial model category.
\end{Lemma}

\Proof
We first remark that
our notion of presentable categories is equivalent to
locally presentable categories in \cite{AR}.
Observe that $\operatorname{S}_{\textup{Tate}}(R)$
is presentable. Since $\Comp(\operatorname{Sh}(\SM,R))$
is presentable
and $\operatorname{S}_{\textup{Tate}}(R)$
can be identified with the functor category
from $\mathfrak{S}'$ to $\operatorname{S}_{\textup{Tate}}(R)$,
thus by \cite[5.5.3.6]{HTT} we see that
$\operatorname{S}_{\textup{Tate}}(R)$ is presentable.
Then according to \cite[3.4.4.2]{HA}
the category $\operatorname{Sp}_{\textup{Tate}}(R)$
of modules over $\operatorname{Sym}(R(1))$
is presentable.
\QED

Let $\operatorname{Comp}(R)$ be the category of chain complexes
of $R$-modules.
There is a combinatorial symmetric monoidal
model structure of $\operatorname{Comp}(R)$
whose weak equivalences are quasi-isomorphisms
and whose fibrations are degreewise surjective maps.
The complex $R$ (concentrated in degree zero) is a cofibrant unit.
This model structure is called the projective model structure (\cite{Ho1}).
There is a symmetric monoidal
functor $\operatorname{Comp}(R)\to \Comp(\operatorname{Sh}(\SM,R))$ which carries a complex $N$ to the constant functor
with value $N$.
For any $A\in \operatorname{Comp}(R)\to \Comp(\operatorname{Sh}(\SM,R))$,
we have the symmetric Tate spectrum
$\{R(1)^{\otimes n}\otimes A\}_{n\in \mathbb{N}}$ such that
$\mathfrak{S}_n$ acts on $R(1)^{\otimes n}\otimes A$
by permutation on $R(1)^{\otimes n}$.
This determines the infinite suspention functor
\[
\Sigma^{\infty}:\Comp(\operatorname{Sh}(\SM,R))\to \operatorname{Sp}_{\textup{Tate}}(R)
\]
which is symmetric monoidal (see \cite[1.4.2.1]{CD2}).
According to \cite[1.2.5, 1.4.2]{CD2},
the composition
\[
\operatorname{Comp}(R)\to \Comp(\operatorname{Sh}(\SM,R))\to \operatorname{Sp}_{\textup{Tate}}(R)
\] is a (symmetric monoidal) left Quillen functor.
By composition,
we also have
\[
L:\SM\to \Comp(\operatorname{Sh}(\SM,R))\stackrel{\Sigma^{\infty}}{\to} \operatorname{Sp}_{\textup{Tate}}(R).
\]

{\it Localizations.}
Now we recall an elegant localization method which transform
model categories into $\infty$-categories (cf. \cite[1.3.4.1, 1.3.1.15, 4.1.3.4]{HA}).
Let $(\mathcal{C},W)$ be a pair of an $\infty$-category
$\mathcal{C}$ and a collection $W$ of edges in $\mathcal{C}$
which contains all degenerate edges.
We say that a map $f:\mathcal{C}\to \mathcal{D}$ exhibits $\mathcal{D}$
as the $\infty$-category obtained from $\mathcal{C}$ by inverting
the edges in $W$ when for any $\infty$-category
$\mathcal{E}$, the functor $f$ induces a fully faithful functor
$\Fun(\mathcal{D},\mathcal{E})\to \Fun(\mathcal{C},\mathcal{E})$
whose essential image consists of functors which sends
edges in $W$ to equivalences in $\mathcal{E}$.
The fibrant replacement $(\mathcal{C},W)\to \mathcal{D}$ of the model category $\textup{Set}_\Delta^+$
of marked simplicial sets (see \cite[3.1]{HTT})
exhibits $\mathcal{D}$ as the $\infty$-category
obtained from $\mathcal{C}$ by inverting
the edges in $W$.
For a model category $\mathbb{M}$, let $\mathbb{M}^c$ be the full subcategory
consisiting of cofibrant objects and $W$ the collection of
edges in $\NNNN(\mathbb{M}^c)$ which correspond to weak equivalences
in $\mathbb{M}^c$.
Then we denote by $\NNNN(\mathbb{M}^c)_{\infty}$ the $\infty$-category
obtained from $\NNNN(\mathbb{M}^c)$ by inverting edges in $W$.
When $\mathbb{M}$ is a combinatorial model category,
$\NNNN(\mathbb{M}^c)_{\infty}$ is a presentable $\infty$-category.
A left Quillen equivalence $\mathbb{M}\to \mathbb{N}$
induces a categorical equivalence
$\NNNN(\mathbb{M}^c)_\infty\to \NNNN(\mathbb{N}^c)_\infty$.
A homotopy (co)limit diagram in $\mathbb{M}$ corresponds to
a (co)limit diagram (see \cite[1.3.4.23, .1.3.4.24]{HA}).
In virtue of \cite[4.1.3.4]{HA}, if $\mathbb{M}$ is a symmetric monoidal model
category, the localization
$\NNNN(\mathbb{M}^c)\to  \NNNN(\mathbb{M}^c)_\infty$ is promoted to
a symmetric monoidal functor $\NNNN(\mathbb{M}^c)^{\otimes}\to \NNNN(\mathbb{M}^c)^\otimes_\infty$ whose underlying functor can be identified with $\NNNN(\mathbb{M}^c)\to  \NNNN(\mathbb{M}^c)_\infty$.
The tensor product $\NNNN(\mathbb{M}^c)_\infty\times \NNNN(\mathbb{M}^c)_\infty\to \NNNN(\mathbb{M}^c)_\infty$ preserves small colimits separately in each variable
since for any $M\in \mathbb{M}^c$,
$(-)\otimes M:\mathbb{M}\to \mathbb{M}$ and $M\otimes (-):\mathbb{M}\to \mathbb{M}$ are left Quillen functors.

Next we apply this localization to the symmetric monoidal left Quillen functor
$\operatorname{Comp}(R)\to \operatorname{Sp}_{\textup{Tate}}(R)$.
Then we have a symmetric monoidal functor of symmetric monoidal presentable $\infty$-categories
\[
\NNNN(\operatorname{Comp}(R)^c)^\otimes_\infty\longrightarrow \NNNN(\operatorname{Sp}_{\textup{Tate}}(R)^c)^\otimes_\infty
\]
which preserves small colimits.
We set
$\mathsf{D}^\otimes(R)=\NNNN(\operatorname{Comp}(R)^c)^\otimes$ and
$\mathsf{Sp}^\otimes_{\textup{Tate}}(R)=\NNNN(\operatorname{Sp}_{\textup{Tate}}(R)^c)^\otimes_\infty$. When we consider the underlying $\infty$-category,
we drop the superscript $\otimes$.
The following Proposition implies that the $\infty$-categories $\mathsf{D}(R)$ and $\mathsf{Sp}_{\textup{Tate}}(R)$ are stable.

\begin{Proposition}
\label{stablemodel}
Let $\mathbb{M}$ be a combinatorial stable model category.
Then the $\infty$-category 
$\NNNN(\mathbb{M}^c)_\infty$ is stable and presentable.
\end{Proposition}

\Proof
The presentability is due to \cite[1.3.4.22]{HA}.

Let $\mathcal{C}=\NNNN(\mathbb{M}^c)_\infty$.
We first observe that
$\mathcal{C}$ is pointed, that is, there is an object which is both initial and final.
According to \cite{D}, the combinatorial model category $\mathbb{M}$
is Quillen equivalent to a combinatorial simplicial model
 category $\mathbb{M}'$.
By \cite[1.3.4.20]{HA} $\calC$ is equivalent to
the nerve $\NNNN((\mathbb{M}')^{\circ})$ where $(\mathbb{M}')^{\circ}$
is the fibrant simplicial category of full subcategory of $\mathbb{M}'$
spanned by cofibrant-fibrant objects. In particular,
the homotopy category of $\calC$
is equivalent to the homotopy category of $\NNNN((\mathbb{M}')^{\circ})$
which is equipped with a structure of a triangulated category.
Let $0$ be a zero object in $\mathbb{M}$ which is cofibrant and fibrant.
We will show that the image $0'$ of $0$ in $\calC$
is a zero object. We prove only that $0'$ is an initial object.
The dual argument shows that $0'$ is also a final object.
By the hammock localization \cite[4.4, 4.7, 5.4]{DK}
together with the equivalence
$\calC \simeq \NNNN((\mathbb{M}')^{\circ})$,
we may identified with $\calC$ with
the nerve of the fibrant replacement of
the hammock localization of $\mathbb{M}^\circ$ (see also \cite[1.3.4.16]{HA}).
Thus for any $X\in \mathbb{M}^\circ$, the homotopy type of the mapping space
from $0$ to $X$ can be calculated by using a simplicial frame of $X$ (cf. \cite[5.4]{Ho1}) and we conclude that the homotopy type is trivial.
Hence $\calC$ is pointed.
Since $\calC$ is presentable,
it has small colimits and limits.
Therefore by \cite[I, 10.12]{DAGn}, it is enough to prove that
the suspension functor $\Sigma$
induces a categorical equivalence $\mathcal{C}\to \mathcal{C}$.
Note that by our assumption and \cite[1.3.4.24]{HA}
the suspention functor induces
an equivalence of the homotopy category
\[
\Sigma: \textup{h}(\calC) \longrightarrow  \textup{h}(\calC).
\]
In particular, $\Sigma:\calC\to \calC$ is essentially surjective.
We claim that $\Sigma:\calC\to \calC$ is fuuly faithful.
It will suffices to show that the suspention functor
induces a homotopy equivalence
$\Map_{\calC}(C,D)\to \Map_{\calC}(\Sigma(C),\Sigma(D))$ for any two objects
$C,D\in \calC$.
Note that $\Map_\calC(C,D)$ is pointed by the zero map
and the natural map
$\Map_\calC(\Sigma (C),D)\to \Omega\Map_\calC(C,D)$ is
a homotopy equivalence.
It follows that
the $n$-th homotopy group $\pi_n(\Map_\calC(C,D))$
can be identified with $\pi_0(\Map_\calC(\Sigma^n (C),D))$.
We conclude that
the map $\pi_n(\Map_\calC(C,D))\to \pi_n(\Map_\calC(\Sigma(C),\Sigma(D)))$
can be identified with the bijective map
$\pi_0(\Map_\calC(\Sigma^n (C),D))\to \pi_0(\Map_\calC(\Sigma^{n+1}(C),\Sigma( D)))$, as desired.
\QED

Let $\mathbf{K}$ be a field of characteristic zero.
Let $\mathbf{HK}$ be the motivic Eilenberg-MacLane spectrum
which is a commutative algebra object in
$\operatorname{Sp}_{\textup{Tate}}(\mathbf{K})$
(see e.g. \cite{RO}).

When $\mathbf{R}$ is a commutative algebra object
in $\operatorname{Sp}_{\textup{Tate}}(\mathbf{K})$
we denote by $\operatorname{Sp}_{\textup{Tate}}(\mathbf{R})$
the category of module objects in $\operatorname{Sp}_{\textup{Tate}}(\mathbf{K})$
over $\mathbf{R}$ (see \cite[Section 4]{ASS}).

According to \cite[1.5.2]{CD2} built on \cite[4.1]{ASS}, there is a combinatorial symmetric monoidal model category structure on $\operatorname{Sp}_{\textup{Tate}}(\mathbf{R})$
such that a morphism is a weak equivalence (resp. fibration) in
$\operatorname{Sp}_{\textup{Tate}}(\mathbf{R})$
if the underlying morphism in $\operatorname{Sp}_{\textup{Tate}}(\mathbf{K})$
is a weak equivalence (resp. fibration).
The base change functor $\operatorname{Sp}_{\textup{Tate}}(\mathbf{K})\to \operatorname{Sp}_{\textup{Tate}}(\mathbf{HK})$ is a symmetric monoidal left Quillen functor.
By inverting by weak equivalences we have a symmetric monoidal functor of symmetric monoidal $\infty$-categories
\[
\mathsf{Sp}^\otimes_{\textup{Tate}}(\mathbf{K})\to \mathsf{Sp}^\otimes_{\textup{Tate}}(\mathbf{HK}):=\NNNN(\textup{Sp}^\otimes_{\textup{Tate}}(\mathbf{HK})^c)^\otimes_\infty
\]
which preserves small colimits.
We remark that $\mathsf{Sp}^\otimes_{\textup{Tate}}(\mathbf{HK})$ is stable
by \cite[4.3.3.17, 8.1.1.4]{HA} and Proposition~\ref{stablemodel}.

\begin{Remark}
There is no reason to assume that
$\mathbf{K}$ is a field of characteristic zero in the above discussion.
We can replace $\mathbf{K}$ by an arbitrary commutative ring $R$.
But in what follows we use the notion of mixed Weil theory which works
over $\mathbf{K}$.
\end{Remark}

\begin{Remark}
\label{motivesconnection}
Let $S$ be the Zariski spectrum of a perfect field $k$.
Let $R$ be an ordinary commutative ring.
Let $\Cor_R$ be the Suslin-Voevodsky's $R$-linear
category of finite correspondences.
Here by an $R$-linear category, we mean
a category enriched over the symmetric monoidal
category of $R$-modules.
An $R$-linear functor means an (obvious) enriched functor.
See \cite{Ke} for the overview of enriched categories.
An object in $\Cor_R$ is a smooth scheme over $S$, that is, an object
in $\SM$. The hom $R$-module
$\Hom_{\Cor_R}(X,Y)$ is a free $R$-module
generated by the set of reduced irreducible closed subscheme
$W\in X\times_k Y$ such that the natural morphism $W\to X$
is finite and its image is an irreducible component of $X$.
The composition
\[
\Hom_{\Cor_R}(X,Y)\otimes_R\Hom_{\Cor_R}(Y,Z)\to \Hom_{\Cor_R}(X,Z),\ \  W\otimes W'\mapsto W' \circ W,
\]
where $W$ and $W'$ are actual reduced irreducible subschemes,
is determined by $W'\circ  W=$ the push-forward by the projection $X\times_k Y\times_k Z\to X\times _k Z$ of the intersection product $(W\times_k Z)\cap (X\times_k W')$.
By the formula $X\otimes Y=X\times_SY$ $\Cor_R$ is a symmetric monoidal
category.
There is a natural map $\SM\to \Cor_{R}$
which sends a smooth scheme $X$ to $X$ and sends morphisms $X\to Y$
to their graphs in $X\times_kY$.
A Nisnevich sheaf of ($R$-modules) with transfers is a contravariant $R$-linear
functor on $\Cor_R$ into the category of $R$-modules, which is a
Nisnevich sheaf on the restriction to $\SM$. Let $\Sh(\Cor_{R})$
be the abelian category of Nisnevich sheaves with transfers.
As the construction of the model category $\textup{Sp}_{\textup{Tate}}(R)$,
in \cite[7.15]{CD1} the symmetric monoidal model
category of $\textup{DM}(S)$ is constructed
(we here employ the notation $\textup{DM}(S)$ in \cite[7.15]{CD1}): we start with the category $\textup{Comp}(\Sh(\Cor_R))$
and take the localization of it by $\mathbb{A}^1$-homotopy equivalence
and stabilize the Tate sphere (this is only the rough strategy, for the detail we refer the reader to \cite{CD1}).
Suppose $R=\mathbf{K}$. There is a left Quillen adjoint symmetric monoidal functor $\textup{Sp}_{\textup{Tate}}(\mathbf{HK})\to \textup{DM}(S)$, which induces
a Quillen equivalence (proved by using alteration \cite[Theorem 68]{RO}, \cite[2.7.9.1]{CD2}).
It gives rise to an equivalence of symmetric monoidal
stable $\infty$-categories
\[
\mathsf{Sp}_{\textup{Tate}}(\mathbf{HK})\to \mathsf{DM}(k):=\textup{N}(\textup{DM}(S)^c)_{\infty}.
\]
Thanks to \cite[2.7.10]{CD2} compact objects and dualizable
objects coincide in $\mathsf{Sp}_{\textup{Tate}}(\mathbf{HK})$.
(We say that an object is dualizable if it have a strong dual
in the sense in loc. cite.)
The full subcategory $\mathsf{Sp}_{\textup{Tate}}(\mathbf{HK})_{\textup{cpt}}$
of the homotopy category of
 $\mathsf{Sp}_{\textup{Tate}}(\mathbf{HK})\simeq \mathsf{DM}(k)$
spanned by compact objects is equivalent to Voevodsky's category
$DM_{gm}(k)$ of geometric motives with coefficients in $\mathbf{K}$.
The triangulated category $DM_{gm}(k)$ is
anti-equivalent to Hanamura's category \cite{HaM}
and
Levine's category \cite{LeM} (with rational coefficients).
\end{Remark}

We summarize the properties of $\mathsf{Sp}_{\textup{Tate}}(\mathbf{HK})\simeq \mathsf{DM}(k)$ as follows:

\begin{Proposition}
\label{goodproperties}
The $\infty$-category
$\mathsf{Sp}_{\textup{Tate}}(\mathbf{HK})\simeq \mathsf{DM}(k)$
is stable and presentable. Moreover,
it is compactly generated (cf. \cite[5.5.7.1]{HTT}). Both compact objects and dualizable obejcts
coincide.
\end{Proposition}

\Proof See Proposition~\ref{stablemodel} and Remark~\ref{motivesconnection}.
\QED

{\it Mixed Weil cohomologies.}
Suppose that the base scheme $S$ is a perfect field $k$.
Let $E$ be a mixed Weil theory in the sense of
Cisinski-D\'eglise \cite[Section 2.1]{CD2}.
A mixed Weil theory is a presheaf $E$ on $\SM$ (or the category of
affine smooth $k$-schemes) of commutative differential graded $\mathbf{K}$-algebras which satisfies $\mathbb{A}^1$-homotopy invariance,
the descent property and axioms on dimension, stability, K\"unneth
formula (see for the detail \cite[2.1.2]{CD2}).
For example, in loc. cite., it is shown that
algebraic and analytic de Rham cohomologies, rigid cohomology,
and $l$-adic \'etale cohomology are mixed Weil theories.
To a mixed Weil theory $E$ we associate a commutative algbera
object $\mathbf{E}$ in $\operatorname{Sp}_{\textup{Tate}}(\mathbf{K})$,
that is, a commutative ring spectrum
(see \cite[2.1.5]{CD2}).
Let $\mathbf{HK}\otimes_{\mathbf{K}}\mathbf{E}$ be the (derived) tensor
product which is a commutative algebra object in $\operatorname{Sp}_{\textup{Tate}}(\mathbf{K})$
(see \cite[2.7.8]{CD2} and its proof).
By \cite[2.7.6]{CD2}, the natural homomorphism
$\mathbf{E}\to \mathbf{HK}\otimes_{\mathbf{K}}\mathbf{E}$
(induced by the structure homomorphism
$\mathbf{K}\to \mathbf{HK}$) is an isomorphism in the homotopy category
of commutative algebra objects.
The homomorphism $\mathbf{E}\to \mathbf{HK}\otimes_{\mathbf{K}}\mathbf{E}$ determines a symmetric monoidal functor $\operatorname{Sp}^\otimes_{\textup{Tate}}(\mathbf{E})\to \operatorname{Sp}^\otimes_{\textup{Tate}}(\mathbf{HK}\otimes_{\mathbf{K}}\mathbf{E})$ which is left Quillen.
The induced symmetric monoidal functor $\rho:\mathsf{Sp}^\otimes_{\textup{Tate}}(\mathbf{E})\to \mathsf{Sp}^\otimes_{\textup{Tate}}(\mathbf{HK}\otimes_{\mathbf{K}}\mathbf{E})$ is an equivalence (since the underlying functor is a categorical
equivalence).
Similarly, there is a symmetric monoidal functor
$\mathsf{Sp}^\otimes_{\textup{Tate}}(\mathbf{HK})\to \mathsf{Sp}^\otimes_{\textup{Tate}}(\mathbf{HK}\otimes_{\mathbf{K}}\mathbf{E})$
determined by the natural homomorphism $\mathbf{HK}\to \mathbf{HK}\otimes_{\mathbf{K}}\mathbf{E}$.
Composing these functors we obtain
\[
\mathsf{D}^\otimes(\mathbf{K})\to \mathsf{Sp}^\otimes_{\textup{Tate}}(\mathbf{K}) \to \mathsf{Sp}^\otimes_{\textup{Tate}}(\mathbf{HK})\simeq \mathsf{DM}^{\otimes}(k) \to \mathsf{Sp}^\otimes_{\textup{Tate}}(\mathbf{HK}\otimes_{\mathbf{K}}\mathbf{E})\stackrel{\rho^{-1}}{\to} \mathsf{Sp}^\otimes_{\textup{Tate}}(\mathbf{E})
\]
where $\rho^{-1}$ is a homotopy inverse of $\rho$.

\begin{Lemma}
\label{stableequiv}
Let $\phi:\mathcal{C}\to \mathcal{D}$ be an exact functor of stable
$\infty$-categories.
Let $\textup{h}(\mathcal{C})$ and $\textup{h}(\mathcal{D})$
be the homotopy categories of $\mathcal{C}$ and $\mathcal{D}$ respectively.
Suppose that $\textup{h}(\phi):\textup{h}(\mathcal{C})\to \textup{h}(\mathcal{D})$
is a categorical equivalence of ordinary categories.
Then $\phi$ is a categorical equivalence.
\end{Lemma}

\Proof
It is clear that $\phi$ is essentially surjective.
It suffices to show that
for $M,N\in \mathcal{C}$,
$\phi$ induces an equivalence 
\[
\Map_{\mathcal{C}}(M,N)\to \Map_{\mathcal{D}}(\phi(M),\phi(N))
\]
in $\SSS$.
We are reduced to proving that the composition
\begin{eqnarray*}
\pi_0(\Map_{\mathcal{C}}(\Sigma^n M, N)))\simeq &\pi_n&(\Map_{\mathcal{C}}(M,N))  \\
&\to& \pi_{n}(\Map_{\mathcal{D}}(\phi(M),\phi(N))) \simeq \pi_0(\Map_{\mathcal{C}}(\Sigma^n \phi(M), \phi(N)))
\end{eqnarray*}
is a bijective where $\pi_n(-)$ denotes the $n$-th homotopy group
and $\Sigma$ is the suspention functor that is compatible with $\phi$.
Now our assertion follows from our assumption.
\QED

\begin{Lemma}
\label{startgoal}
The composition $\mathsf{D}^\otimes(\mathbf{K}) \to \mathsf{Sp}^\otimes_{\textup{Tate}}(\mathbf{E})$ is an equivalence of symmetric monoidal $\infty$-categories.
\end{Lemma}

\Proof
It is enough to show that the underlying functor is a categorical equivalence.
By Lemma~\ref{stableequiv}
it suffices to prove that
the induced functor of homotopy
categories
$\textup{h}(\mathsf{D}(\mathbf{K}))\to \textup{h}(\mathsf{Sp}_{\textup{Tate}}(\mathbf{E}))$ is an equivalence.
The right adjoint of this functor is described as
$D_{\mathbf{A}^1}(k,\mathbf{E})=\textup{h}(\mathsf{Sp}_{\textup{Tate}}(\mathbf{E})) \to D(\mathbf{K})=\textup{h}(\mathsf{D}(\mathbf{K}))$ given by
$M\mapsto \mathbf{R}\Hom_{\mathbf{E}}(\mathbf{E},M)$ where 
we use the notation $D_{\mathbf{A}^1}(k,\mathbf{E})$,
$D(\mathbf{K})$ and $\mathbf{R}\textup{Hom}_{\mathbf{E}}(\mathbf{E},M)$
in \cite{CD2} (namely, the right adjoint is given by the ``Hom complex''
$\mathbf{R}\textup{Hom}_{\mathbf{E}}(\mathbf{E},M)$ in $\textup{h}(\mathsf{Sp}_{\textup{Tate}}(\mathbf{E}))$).
This right adjoint is an equivalence by \cite[2.7.11]{CD2} and thus
$\textup{h}(\mathsf{D}(\mathbf{K}))\to \textup{h}(\mathsf{Sp}_{\textup{Tate}}(\mathbf{E}))$ is so.
\QED

Let $H\mathbf{K}$ be the (not motivic)
Eilenberg-MacLane commutative ring spectrum of $\mathbf{K}$
in $\SP$.

\begin{Proposition}
\label{monoidalequivalencemodule}
There is an equivalence
$\Mod_{H\mathbf{K}}^\otimes \to  \mathsf{D}^\otimes (\mathbf{K})$
of symmetric monoidal $\infty$-categories.
\end{Proposition}

\Proof
This immediately follows from \cite[8.1.2.13]{HA}.
\QED

\begin{Remark}
There is no need to assume that $\mathbf{K}$ is a field.
The proof is valid for any commutative ring.
\end{Remark}

\begin{Definition}
By Proposition~\ref{monoidalequivalencemodule} and
Lemma~\ref{startgoal},
we obtain a symmetric monoidal functor
\[
\mathsf{R}_E:\mathsf{Sp}_{\textup{Tate}}^\otimes(\mathbf{HK})\simeq \mathsf{DM}^\otimes(k)
\to \mathsf{Sp}^\otimes_{\textup{Tate}}(\mathbf{HK}\otimes_{\mathbf{K}}\mathbf{E})\stackrel{\rho^{-1}}{\to} \mathsf{Sp}^\otimes_{\textup{Tate}}(\mathbf{E})\simeq \Mod_{H\mathbf{K}}^\otimes.
\]
We refer to is as the realization functor assosiated to $E$.
\end{Definition}

\subsection{The construction of motivic Galois groups}
For a mixed Weil thoery $E$, we have
\[
\xymatrix{
\mathsf{Sp}^\otimes_{\textup{Tate}}(\mathbf{HK})\simeq \mathsf{DM}^\otimes(k) \ar[r]^(0.6){\mathsf{R}_E} & \Mod_{H\mathbf{K}}^\otimes.
}
\]
For example, suppose that $S=\Spec k$ is the Zariski spectrum
of a field of characteristic
zero and $\mathbf{K}=k$.
Let $E$ be the mixed Weil theory of algebraic de Rham cohomology
and $L(X)$ the image of smooth scheme $X\in \SM$ in $\mathsf{Sp}^\otimes_{\textup{Tate}}(\mathbf{HK})$.
Then $\mathsf{R}_E$ carries $L(X)$ to the dual of the complex computing the de Rham cohomology of $X$.

By Remark~\ref{motivesconnection}, in
$\Mod^\otimes_{H\mathbf{K}}$ and $\mathsf{Sp}^\otimes_{\textup{Tate}}(\mathbf{HK})$,
compact objects and dualizable objects coincide respectively.
This diagram induces the diagram of full subcategories of
dualiziable objects whose underlying $\infty$-categories are small stable idempotent complete
$\infty$-categoires ($\heartsuit$):
\[
\xymatrix{
\mathsf{Sp}^\otimes_{\textup{Tate}}(\mathbf{HK})_{\vee}\simeq \mathsf{DM}_{\vee}^\otimes(k) \ar[r]^(0.6){\mathsf{R}_E} & \PMod_{H\mathbf{K}}^\otimes
}
\]
where $\mathsf{R}_E$ is the restriction of the realization functor
(we abuse notation).

\begin{Definition}
\label{mgg}
We apply Theorem~\ref{supermain} to the realization functor
$(\heartsuit)$ and obtain a derived affine group scheme $\mathsf{MG}_{E}$ over $H\mathbf{K}$
which we shall call the derived motivic Galois group associated to the mixed Weil theory $E$.
There is a diagram of symmetic monoidal stable idempotent complete $\infty$-categories
\[
\xymatrix{
\DM_\vee^\otimes(k) \ar[dr]_{\mathsf{R}_E} \ar[rr] &  & \PRep_{\mathsf{MG}_E}^\otimes \ar[dl] \\
  & \PMod^\otimes_{H\mathbf{K}} & 
}
\]
where $\PRep_{\mathsf{MG}_E}^\otimes$ is the symmetric monoidal stable idempotent complete $\infty$-category of perfect representations of $\mathsf{MG}_{E}$ (see Appendix A.6) and $\PMod_{\mathsf{MG}_E}^\otimes \to \PMod^\otimes_{H\mathbf{K}}$ is the forgetful functor.
When $E$ is clear, we often write $\mathsf{MG}$ for $\mathsf{MG}_E$.
If we let $\mathsf{MG}_E=\Spec B_E$, then we can choose $B_E$ to be
a commutative differential graded $\mathbf{K}$-algebra
$B_E$
by virtue of the well-known categorical equivalence
between the $\infty$-category of commutative $H\mathbf{K}$-ring spectra
and that of commutative differential graded $\mathbf{K}$-algebras
(cf. e.g. \cite[8.1.4.11]{HA}).
\end{Definition}

\begin{Theorem}
\label{derivedgroup}
The derived affine group scheme $\mathsf{MG}_E=\Spec B_E$
has the universality described in Theorem~\ref{supermain}
and represents
the automorphism group functor $\Aut(\mathsf{R}_E)$.
\end{Theorem}

\begin{Remark}
Since $\mathbf{K}$ is a field of characteristic zero,
to work with $\mathsf{MG}_E$,
we may employ complicial algebraic geometry \cite[II, 2.3]{HAG2}.
But when one wants to apply our tannakization
to the integral Betti realization
and obtain motivic Galois group over $H\mathbb{Z}$,
we need the brave new derived algebraic geometry \cite[II, 2.4]{HAG2}, \cite{DAGn}.
\end{Remark}

{\it Variants.}
Theorem~\ref{supermain} is quite powerful.
We can also construct a derived affine group scheme
from any symmetric monoidal (full) subcategory in $\DM_{\vee}^\otimes(k)$.
Let $\mathsf{S}^\otimes\subset\DM_\vee^\otimes(k)$ be a symmetric monoidal
full subcategory.
In virtue of Theorem~\ref{supermain}
the composite
\[
\mathsf{S}^\otimes \hookrightarrow \DM_\vee^\otimes(k)\stackrel{\mathsf{R}_E}{\longrightarrow} \PMod^\otimes_{H\mathbf{K}}
\] 
yields a derived affine group scheme $\mathsf{MG}_E(\mathsf{S}^\otimes)$
over $H\mathbf{K}$.
Full subcategories of mixed Tate motives, Artin motives and so on have been
very important examples.
As mentioned in Introduction, we will investigate
the tannakizations of these full subcategories in a
separate paper \cite{Bar}.

Let $X$ be a smooth scheme over $k$. Let $m$ be an integer.
Let $\DM_\vee^\otimes(k)_{X(m)}$ denotes the smallest
symmetric monoidal idempotent complete stable subcategory
which contains $L(X)(m)$.
The underlying stable $\infty$-category
is the smallest stable subcategory which contains
$(L(X)(m))^{\otimes n}$ ($n \ge 0$)
and is closed under retracts.
In this case, we write $\mathsf{MG}_E(X(m))$
for $\mathsf{MG}_E(\DM_\vee^\otimes(k)_{X(m)})$.

\begin{Proposition}
There exists a natural equivalence of derived affine group schemes
\[
\mathsf{MG}_E\stackrel{\sim}{\longrightarrow} \lim_{(X,m)}\mathsf{MG}_E(X(m))
\]
where the right-hand side is the (small) limit of derived affine
group schemes, and pairs
$(X,m)$ run over smooth projective schemes $X$ and integers $m\in \ZZ$.
\end{Proposition}

\Proof
It is enough to show that the colimit
\[
\textup{colim}_{(X,m)}\DM_\vee^\otimes(k)_{X(m)}
\]
in $\CAlg(\uCat)$
is equivalent to $\DM_\vee^\otimes(k)$.
It will suffice to prove that
the colimit of the diagram of
$\{\DM_\vee(k)_{X(m)}\}_{(X,m)}$
in $\textup{Cat}_\infty$ is equivalent to the $\infty$-category
$\DM_\vee^\otimes(k)$.
We are reduced to showing that for any $M\in \DM_\vee^\otimes(k)$
there exists $L(X)(m)$ such that
$M$ belongs to $\DM_\vee^\otimes(k)_{X(m)}$.
Note that by Proposition~\ref{goodproperties}
there exists a finite collection of
objects $\{L(X_1)(m_1),\ldots,L(X_r)(m_r)\}$
such that $M$
lies in the smallest
stable subcategory
which contains all $L(X_i)(m_i)$ and is closed under retracts.
Therefore we easily see that there exists $L(X)(m)$ such that
$M\in \DM_\vee^\otimes(k)_{X(m)}$
\QED

{\it Truncated affine group schemes.}
We can obtain an ordinary affine group scheme over $\mathbf{K}$
from $\mathsf{MG}_E$ or variants.
For this purpose, let $\mathsf{dga}_\mathbf{K}$
be the category of commutative differential graded $\mathbf{K}$-algebras.
By virtue of \cite[2.2.1]{Hi} or \cite[8.1.4.10]{HA}, there is a combinatorial
model category structure
on $\mathsf{dga}_\mathbf{K}$, in which
weak equivalences are quasi-isomorphisms (of underlying complexes), and fibrations
are those maps which induce levelwise surjective maps.
Let $\mathsf{dga}^{\ge 0}_\mathbf{K}$
be the full subcategory of $\mathsf{dga}_\mathbf{K}$
spanned by those objects such that $A^i=0$ for any $i<0$
(here we use the cohomological indexing).
According to \cite[2.2.1]{Hi}, there is a
combinatorial model category structure
on $\mathsf{dga}^{\ge0}_\mathbf{K}$, in which
weak equivalences are quasi-isomorphisms, and fibrations
are those maps which induce levelwise surjective maps (here
one can choose
a Sullivan algebra as a cofibrant replacement, see e.g. \cite{Hess}).
It gives rise to a Quillen adjunction
\[
\tau:\mathsf{dga}_\mathbf{K}\rightleftarrows \mathsf{dga}^{\ge0}_\mathbf{K}:\iota
\]
where the right adjoint $\iota:\mathsf{dga}^{\ge0}_\mathbf{K}\to \mathsf{dga}_\mathbf{K}$ is the inclusion. The left adjoint sends $A$ to its quotient
by differential graded ideal generated by elements $a\in A^i$ for $i<0$.
By the localization,
we have the induced adjunction
\[
\CAlg_{H\mathbf{K}}\simeq \NNNN(\mathsf{dga}^c_\mathbf{K})_\infty\rightleftarrows \NNNN((\mathsf{dga}^{\ge0}_\mathbf{K})^c)_\infty
\]
(cf. \cite[1.3.4.26]{HA} and adjoint functor theorem \cite[5.5.2.9]{HA}), where the right adjoint is fully faithful. For the first equivalence,
see \cite[8.1.4.11]{HA}.
Hence by adjunction, $\Spec \tau B_E$
represents $\NNNN((\mathsf{dga}^{\ge0}_\mathbf{K})^c)_\infty\hookrightarrow \CAlg_{H\mathbf{K}}\stackrel{\Aut(\mathsf{R}_E)}{\to}\Grp(\mathcal{S})$
(recall $\mathsf{MG}_E=\Spec B_E$).
Next let $\mathsf{dga}^{0}_\mathbf{K}$ be the full subcategory
of $\mathsf{dga}_\mathbf{K}$ spanned by objects such that $A^i=0$ for
$i\ne 0$, that is, the category of commutative $\mathbf{K}$-algebras.
Then there is a natural adjunction
$\NNNN(\mathsf{dga}^0_\mathbf{K}) \rightleftarrows \NNNN((\mathsf{dga}^{\ge0}_\mathbf{K})^c)_\infty$,
where the right adjoint carries $A$ to the (homologically)
connective cover of $A$, i.e.
$H^0(A)$, and the left adjoint is the natural inclusion.
Note that since $\mathbf{K}$ is a field there is a natural isomorphism
$H^0(A\otimes B)\simeq H^0(A)\otimes H^0(B)$ for $A,B\in \mathsf{dga}^{\ge0}_\mathbf{K}$.
Therefore, the affine scheme
$\Spec H^0(\tau B_E)$ inherits a group structure
from $\Spec\tau B_E$.
We refer to the affine group scheme (i.e. pro-algebraic group)
\[
MG_E=\Spec H^0(\tau B_E)
\]
as the underived motivic Galois group.
The affine group scheme $MG_E$ has the following important property:

\begin{Theorem}
\label{ordinarygroup}
Let $K$ be a $\mathbf{K}$-field.
Let $\overline{\Aut(\mathsf{R}_E)}(K)$ be the group of
isomorphism classes of automorphisms of $\mathsf{R}_E$,
that is, $\pi_0(\Aut(\mathsf{R}_E)(HK))$.
Then there is a natural isomorphism of groups
\[
MG_E(K) \simeq \overline{\Aut(\mathsf{R}_E)}(K)
\]
where $MG_E(K)$ denotes the group of $K$-valued points.
These isomorphisms are functorial among $\mathbf{K}$-fields
in the obvious way.
\end{Theorem}

\Proof
We first suppose that $K=\mathbf{K}$. Let $B_E'=\tau B_E$.
There is a natural morphism
$\Spec B'_E \to \Spec H^0(B'_E)$
corresponding to $H^0(B'_E)\to B'_E$.
Here $\Spec(-)$ is considered to be a functor
$\NNNN((\mathsf{dga}_\mathbf{K}^{\ge 0})^c)_\infty\to \Grp(\mathcal{S})$.
Let $u:\Spec \mathbf{K}\to \Spec H^0(B'_E)$ be the unit morphism
and let $\Spec C=\Spec \mathbf{K}\times_{\Spec H^0(B'_E)}\Spec B'_E$
be the associated (homotopy) fiber in $\NNNN((\mathsf{dga}^{\ge0}_\mathbf{K})^c)_\infty^{op}$.
Observe that $H^n(B'_E)$ is a free $H^0(A)$-module for any $n \ge0$.
To see this, note that
the commutative Hopf graded algebra $H^*(B'E)$
induces a coaction $H^n(B'_E)\to H^n(B'_E)\otimes_{\mathbf{K}}H^0(B'_E)$
of the commutative Hopf algebra $H^0(B'_E)$.
This action commutes with the coaction of $H^0(B'_E)$ on itself.
Consequently, the quasi-coherent module $H^n(B'_E)$ on
$\Spec H^0(B'_E)$ descends to $\Spec \mathbf{K}\simeq \Spec H^0(B'_E)/\Spec H^0(B'_E)$. It follows that $H^n(B'_E)$ has the form $L\otimes_{\mathbf{K}}H^0(B'_E)$ where $L$ is a vector space.
The freeness implies that
$H^n(C)\simeq \mathbf{K}\otimes_{H^0(B'_E)}H^n(B'_E)$.
In particular, $H^0(C)\simeq \mathbf{K}$.
Hence by \cite[VIII, 4.4.8]{DAGn},
for any usual $\mathbf{K}$-algebra $R$, $\Map_{\NNNN((\mathsf{dga}^{\ge0}_\mathbf{K})^c)_\infty^{op}}(\Spec R,\Spec C)$ is connected.
Thus the full subcategory of $\Map_{\NNNN((\mathsf{dga}^{\ge0}_\mathbf{K})^c)_\infty^{op}}(\Spec R,\Spec B'_E)$ spanned by morphisms $\Spec R\to \Spec B'_E$
lying over $u$, is connected.
If we replace $u$ by another $\mathbf{K}$-valued point $v$ of $MG_E$
via a translation by group action, the same conclusion holds.
Therefore, we have a natural isomorphism $\overline{\Aut(\mathsf{R}_E)}(\mathbf{K})\simeq MG_E(\mathbf{K})$.
For a general $\mathbf{K}$-field $K$, if we replace $B'_E$ by the base change
$B'_E\otimes_{\mathbf{K}}K$, then the same argument works.
\QED

\begin{Remark}
\label{abelianmotives}
The tannakian view of motives is originated from Grothendieck's idea.
For the original ideas of motivic Galois groups and motivations,
we refer the reader to \cite{And}.
The guiding principle behind our work is that
the stable $\infty$-category of mixed motives (or so-called geometric motives)
should naturally constitute a ``tannakian category''
in the setting of $\infty$-categories.
It is considered to be a version of the original idea, which is generalized to
the realm of higher category theory.
Arguably, a conjectural abelian (furthermore tannakian) category
of mixed motives is defined as 
the heart of $\mathsf{DM}_\vee^\otimes(k)$
endowed with a conjectural motivic $t$-structure.
Here a motivic $t$-structure is a nongenerate $t$-structure
on the homotopy category of $\mathsf{DM}_\vee(k)$,
such that $\otimes:\mathsf{DM}_\vee(k)\times \mathsf{DM}_\vee(k)\to \mathsf{DM}_\vee(k)$ and the realization functor are $t$-exact.
The existence of a motivic $t$-structure is a hard problem,
and recently Beilinson \cite{Be} shows that the existence
of a motivic $t$-structure implies Grothendieck's standard conjectures
(cf. \cite[Chapitre 5]{And})
in characteristic zero.
Conversely, Hanamura \cite{Ha3} proves that
a ``generalized standard conjectures'' including Beilinson-Soul\'e
vanishing conjecture imply the existence of a
motivic $t$-structure.
(It is worth remarking that a construction of
a motivic Galois group for numerical pure motives also needs
the standard conjectures, see \cite{And}.)
These conjectures in full generality are largely inaccessible
by now.
The idea of us is to start with the $\infty$-category
$\mathsf{DM}^\otimes(k)$ endowed with the realization functor
into $\infty$-category of complexes,
partly motivated by ``derived tannakian philosophy''.
The reader might raise an objection to our construction of the motivic
Galois group as a derived affine group scheme.
(But we can extract a usual group scheme from it as above.)
We do not think that this is the drawback.
Rather, the derived affine group scheme $\mathsf{MG}=\Spec B$ should
capture the interesting new data of ``highly structured'' category
$\mathsf{DM}(k)$ of mixed motives which may not
arise from a conjectural abelian category of mixed motives.
Suppose that a motivic $t$-structure exists
and let $\mathcal{MM}$ be its heart.
Let $D^b(\mathcal{MM})$ be the bounded derived category (if exists)
and let $D^b(\mathcal{MM})\to DM_{gm}(k)$ be the natural functor.
The problem whether or not $D^b(\mathcal{MM})\to DM_{gm}(k)$ is an equivalence
is mysterious.
Thus, at least a priori, we can think that
$\mathsf{DM}^\otimes_{\vee}(k)$
has richer information than $\mathcal{MM}$.
We morally think of the part of higher and lower homotopy data of $\mathsf{MG}$
as the data of $\mathsf{DM}(k)$
which can not be determined by the abelian category $\mathcal{MM}$.
Beside, this might reveal new insights on the motivic Galois
group of a conjectural abelian category of mixed motives.

In the case of mixed Tate motives,
Beilinson-Soul\'e vanishing conjecture implies the existence of a motivic
$t$-structure on the triangulated subcategory of mixed Tate motives,
by the work of Kriz-May
\cite{KM}, Levine \cite{Lev}.
In \cite{BK} and \cite{KM}, the bar construction
of a motivic dg-algebra is used, and it yields a derived affine group scheme.
Recently, using bar constructions Spitzweck
has constructed the derived affine group scheme
such that its representation category is equivalent to
the ($\infty$-)category of (integral) mixed Tate motives, see \cite{Spi}.
This construction can be viewed as Beilinson-Soul\'e vanishing conjecture-free 
and $K(\pi,1)$-property-free approach.
In \cite{Bar}, as mentioned before, we study the tannakzaiton of $\infty$-category of
mixed Tate motives, which is related to the so-called 
motivic Galois group for mixed Tate motives.
\end{Remark}

\begin{Remark}
There is the natural functor
$\DM_{\vee}^\otimes(k) \longrightarrow \PRep^\otimes_{\mathsf{MG}_E}$.
It seems reasonable to conjecture that
this functor is an equivalence.
This conjecture is a refinement of \cite[22.1.4.1 (ii)]{And}
which says that
the realization functor is conservative, that is, $\mathsf{R}_E(M)=0$
implies that $M=0$.
Of course, this functor is universal among functors into
the  $\infty$-categories of complexes
of the representations of affine groups over $\mathbf{K}$.
Namely, let $f:\mathsf{DM}^\otimes_{\vee}(k)\to \PRep^\otimes_{G}$ be a functor which commutes with functors
to $\PMod^\otimes_{H\mathbf{K}}$ where $G$ is
a usual affine group scheme over $\mathbf{K}$ (considered as the derived
affine group scheme). Then there exists a homomorphism
$G\to \mathsf{MG}_E$ which induces
$\PMod^\otimes_{\mathsf{MG}_E}\to \PMod^\otimes_G$ such that
the composition 
$\mathsf{DM}_{\vee}^\otimes(k) \to \PRep^\otimes_{\mathsf{MG}_E}\to \PRep^\otimes_G$ is equivalent to $f$.
An example of such $G$ we should keep in mind is the Tannaka dual of
the abelian category of finite dimensional
continuous $l$-adic representations
of the absolute Galois group (when $\mathbf{K}=\mathbb{Q}_l$ and $E$ is
the mixed Weil theory of $l$-adic \'etale cohomology).
Another important example is the Tannaka dual of the abelian category of
mixed Hodge structures.

\begin{Remark}
There has been Nori's abelian category of mixed motives (see \cite{Ara})
and its motivic Galois group $\mathsf{MG}_{\textup{Nori}}$.
It is natural to consider
that the relationship between our $\mathsf{MG}$ and $\mathsf{MG}_{\textup{Nori}}$.
Our $\mathsf{MG}$ is directly related with $\mathsf{DM}^\otimes(k)$ and the
realization functor, and
this question depends on the relation between 
$\mathsf{DM}^\otimes(k)$
and ($\infty$-categorical setup of) the derived category
of Nori's abelian category, which seems out of reach at the present time.
\end{Remark}

\end{Remark}

\section{Other examples}

In this Section,
we will present some other examples for the applications of tannakizations.
To avoid getting this Section long, we only mention examples
 which one can define quickly.

\subsection{Perfect complexes of derived stacks}

Let $R$ be a commutative ring spectrum.
Let $\XXX$ be a derived stack over $R$ (for this notion, we refer to
\cite{HAG2}, \cite{DAGn}, or \cite{Bar}).
Let $\Perf^\otimes(\XXX)$ be
the symmetric monoidal stable $\infty$-category of
perfect complexes on $\XXX$.
Here we define $\Perf^\otimes(\XXX)$ to be
the limit $\lim_{\Spec A\to \XXX}\PMod^\otimes_A$ in $\CAlg(\textup{Cat}_\infty)$
where $\Spec A\to \XXX$ run over smooth morphisms with $A\in \CAlg_R$.
Let $p:\Spec R\to \XXX$ be a morphism of derived stacks over $R$.
We have the pullback functor
\[
p^*:\Perf^\otimes(\XXX)\to \PMod_{R}^\otimes \simeq \Perf^\otimes(\Spec R)
\]
which is an $R$-linear symmetric monoidal exact functor.
It gives rise to its tannakziation; a derived affine group scheme over $R$.
We can think this as a generalization of
bar constructions of commutative ring spectra.
In \cite{Bar} we study this issue in detail.

\subsection{Topological spaces}
Let $R$ be a connective commutative ring spectrum.
Let $S$ be a topological space
which we regard as an object in $\mathcal{S}$.
Let $p:\Delta^0\to S$ denote a point.
We can view $S$ as a constant functor belonging to
$\Fun(\CAlg^{\textup{con}}_R,\mathcal{S})$.
Let $\operatorname{Perf}^{\otimes}(S)$ be 
the limit
$\lim_{\Spec R\to S} \PMod^{\otimes}_R$
where
$\Spec R\to S$ run over $\Mod_{\Fun(\CAlg^{\textup{con}}_R,\mathcal{S})}(\Spec R,S)$.
We may think of $\operatorname{Perf}^{\otimes}(S)$ as the
symmetric monoidal $\infty$-category of perfect complexes on $S$ with $R$-coefficients.
The symmetric monoidal $\infty$-category
$\operatorname{Perf}^{\otimes}(S)$ is a small stable idempotent complete
$\infty$-category.
Then the prescribed point $p:\Delta^0\to S$ induces
\[
\xymatrix{
\Perf^\otimes(S) \ar[r] & \Perf^\otimes(\Delta^0)
}
\]
where $\Perf^\otimes(\Delta^0)\simeq \Perf^\otimes(\Spec R)\simeq \PMod^\otimes_R$.
We then apply the tannakization functor to this diagram.
We denote by $G(S,p)$ the associated derived affine group scheme
over $R$.

When $R=H\mathbb{Q}$, it would be interesting to compare
the rational homotopy theory and $G(S,p)$ over $H\mathbb{Q}$.
We speculate on the relation to the de Rham homotopy theory.
For simplicity, $S$ is simply connected of finite type.
Let $A_{PL}(S)$ be the polynomial de Rham algebra of $S$
over $\mathbb{Q}$
(see e.g. \cite{BG}).
It is a commutative differential graded $\mathbb{Q}$-algbera.
Since the coefficient is $\mathbb{Q}$,
we may regard $A_{PL}(S)$ as a coconnective commutative ring
spectrum over $H\mathbb{Q}$ (that is, $\pi_i(A_{PL}(S))=0$ for $i>0$).
Let $\operatorname{Spec}(A_{PL}(S))$ be
the functor $\CAlg_{H\mathbb{Q}} \to \mathcal{S}$
corepresentable by $A_{PL}(S)$.
The restriction of $\operatorname{Spec}(A_{PL}(S))$
to $\CAlg_{H\mathbb{Q}}^{\textup{con}}$ is 
a schematization of $S$ (see \cite[VIII, 4.4.2]{DAGn}, \cite{AC}).
There is a natural base point $\rho:\operatorname{Spec}(H\mathbb{Q})\to \operatorname{Spec}(A_{PL}(S))$ induced by $\Delta^0\to S$.
The associated \v{C}ech nerve of $\rho$ determines a simplicial
diagram $\NNNN(\Delta)^{op}\to \Aff_{H\mathbb{Q}}$
which is a derived affine group scheme $G_{PL}(S)$.
Then the relationship with de Rham homotopy theory
should be described by an equivalence $G(S,p)\simeq G_{PL}(S)$
of derived group schemes over $H\mathbb{Q}$
($G_{PL}(S)$ is obtained by the tannakization of the forgetful functor $\PMod_{A_{PL}(S)}^\otimes \to \PMod^\otimes_{H\mathbb{Q}}$).
We hope that
our construction brings a new conceptual insight to rational homotopy theory
and wish to return this issue in the future.

\appendix

\renewcommand{\theTheorem}{A.\arabic{Theorem}}

\section{Derived group schemes.}

\subsection{Derived schemes}
Before proceeding to derived (affine) group schemes,
let us review derived schemes and fix our convention.
Let $R$ be a commutative ring spectrum.
Recall that $\CAlg$ denotes the $\infty$-category of commutative
ring spectra (commutaive algebra objects in
$\SP$, i.e. $E_{\infty}$-rings in \cite{HA}).
We will fix our convention on derived schemes.

Let us recall the notion of \'etale and flat morphisms in $\CAlg$.
We say that a morphism $A\to B$ in $\CAlg$ is
\'etale (resp. flat)
if it has the following properties:
\begin{enumerate}
\item $\pi_0(A)\to \pi_0(B)$ is \'etale (resp. flat),

\item the isomorphism $\pi_0(B)\otimes_{\pi_0(A)}\pi_n(A)\simeq \pi_n(B)$ of abelian groups for any $n\in \ZZ$.

\end{enumerate}

If an \'etale (resp. flat) morphism $A\to B$ induces
a surjective morphism $\Spec \pi_0(B)\to \Spec \pi_0(A)$,
we say that $A\to B$ is \'etale (resp. flat) surjective.

Let $A\to B^\bullet$ be an coaugmented cosimplicial objects in $\CAlg_R$.
We say that $A\to B^\bullet$ is an \'etale hypercover
if for any $n\ge0$, the natural
morphism $\textup{cosk}_{n-1}(B^\bullet)_n\to B^n$
is \'etale surjective, and $A\to B^0$ is \'etale surjective.
Here we abute notation by writing $\textup{cosk}_{n-1}(B^\bullet)_n$
for the coskeleton when we consider $B^\bullet$ 
to be the simplicial object in $(\CAlg_R)^{op}$.

We say that a functor (or presheaf)
$P:\CAlg_{R}\to \SSS$ is a (hypercomplete \'etale)
sheaf if the following two properties hold:
\begin{itemize}
\item if $\{A_{\lambda}\}$ is a finite family of objects
in $\CAlg_{R}$, then $P(\sqcap_{\lambda} A_{\lambda})\simeq \sqcap_{\lambda}P(A_{\lambda})$

\item 
Let $A\to B^\bullet$ be an \'etale hypercover.
Then we have $P(A)\simeq \lim(P(B^{\bullet}))$,
where $\lim(P(B^{\bullet}))$ denotes a limit of the cosimplicial diagram.

\end{itemize}

Let $\Sh(\CAlg_{R}^{et})$ be the full subcategory of
$\Fun (\CAlg_{R},\widehat{\SSS})$ spanned by sheaves.
($\widehat{\SSS}$ is the $\infty$-category of spaces in an enlarged universe.)
For any $A$ in $\CAlg_R$, we define $\Spec A$
to be a functor $\CAlg_R\to \mathcal{S}$ corepresentable by $A$.
This functor is a sheaf. Namely, $\Spec A$ belongs to $\Sh(\CAlg_{R}^{et})$.
We shall refer to $\Spec A$ as the derived affine scheme (over $R$) associated to $A$.
Let $\Aff_{R}\subset \Sh(\CAlg_{R}^{et})$ be the full subcategory spanned by derived affine schemes over $R$.
Yoneda's Lemma implies that $\Aff_{R}\simeq (\CAlg_{R})^{op}$.
If $R$ is the sphere spectrum, then we usually write $\Aff$ for $\Aff_{R}$.

A derived scheme is informally a ``geometric object'' which is
``Zariski locally'' isomorphic to a derived affine scheme.
In \cite{DAGn}, Lurie develops the approach of ringed $\infty$-topoi
to the definition of derived schemes and derived Deligne-Mumford stacks.
We here take the definition of derived schemes which is similar to
To\"en-Vezzosi \cite{HAG2}.
A derived scheme 
over $R$ which has affine diagonal
is a sheaf (that is, a contravariant functor which satisfies the descent condition as above) $X:\Aff_R^{op}\to \widehat{\mathcal{S}}$
which  has  the following properties (i) and (ii), 
\begin{enumerate}
\renewcommand{\labelenumi}{(\roman{enumi})}

\item for any two morphisms (natural transformations)
$a:\Spec A \to X$ and $b:\Spec B\to X$
with derived affine schemes $\Spec A$ and $\Spec B$ over $R$,
then the fiber product $\Spec A\times_{X}\Spec B$ is
representable by a derived affine scheme $\Spec C$,

\item there exist the disjoint union of derived affine
schemes $\sqcup_{\lambda\in I}\Spec A_{\lambda}$ and 
a morphism $p:\sqcup_{\lambda\in I}\Spec A_{\lambda}\to X$
such that for any $q:\Spec B\to X$ and any $\lambda\in I$,
the base change
$\sqcup_{\lambda}\Spec C_{\lambda} \to \Spec B$
is an \'etale morphism
and
it induces an open immersion
$\Spec \pi_0(C_{\lambda})\to \Spec \pi_0(B)$ for each $\lambda\in I$,
and a surjective morphism
$\sqcup_{\lambda}\Spec \pi_0(C_{\lambda}) \to \Spec \pi_0(B)$ of ordinary
schemes,
where $\Spec C_\lambda:=\Spec A_{\lambda}\times_{X}\Spec B$.

\end{enumerate}

We denote by $\operatorname{Sch}_R$ the full subcategory spanned by derived
schemes over $R$. (We assume that
all derived schemes have affine diagonal.)

We shall refer to \cite[II, 2,4]{HAG2}, \cite{DAGn} for the generalities on derived
schemes and derived stacks.

\begin{Remark}
In this paper we work with the derived 
algebraic geometry over nonconnective commutative ring spectra
(this point is relevant to motivic applications).
\end{Remark}

\subsection{Derived group schemes}
A (ordinary) group scheme over a scheme $S$
is a scheme $G$ which is endowed with morphisms
$S\to G$ and $G\times_S G\to G$ that satisfies the usual group axioms.
If one employs the functorial point of view,
then a group scheme is a group-valued functor on the category of commutative rings, which is representable by a scheme.
The notion of derived group schemes is similar to that of
group schemes.
The point is that to define the notion of
derived group schemes we will replace the ordinary category of
commutative rings by $\CAlg$.
As the case of derived schemes,
the notion of group-valued functors on $\CAlg$ is not useless.
We should treat functors into group objects in $\mathcal{S}$.
We first recall the notion of group objects in $\infty$-categorical settings
(these are also commonly called group-like $A_{\infty}$-spaces in operadic contexts).
We refer to \cite{TH} \cite{Spi} for accounts of this subject including
related notions.

\begin{Definition}
\label{monoidgroup}
Let $\mathcal{C}$ be an $\infty$-category which admits finite limits.
A monoid object is a map $f:\NNNN(\Delta)^{op}\to \mathcal{C}$
having the property: $f([0])$ is a final object,
and for each $n\in \mathbb{N}$, 
inclusions $\{i-1,i\}\hookrightarrow[n]$ for $1\le i\le n$
induce an equivalence
\[
f([n])\to f([1])\times \ldots \times f([1])
\]
where the right hand side is the $n$-fold product.
We denote by $\textup{Mon}(\mathcal{C})$ the full subcategory of
$\Fun(\NNNN(\Delta)^{op},\mathcal{C})$ spanned by monoid objects.

A groupoid object in $\mathcal{C}$ is a functor
$f:\textup{N}(\Delta)^{op}\to \mathcal{C}$ with the following property:
for every $n$ and every partition $[n]=S\cup S'$ such that $S\cap S'$ has
one element which we denote by $s$, the diagram
\[
\xymatrix{
f([n]) \ar[r] \ar[d] & f(S) \ar[d] \\
f(S') \ar[r] & f(\{s\})
}
\]
is a pullback diagram in $\mathcal{C}$ (see \cite[6.1.2]{HTT}).
We say that a groupoid object $f:\textup{N}(\Delta)^{op}\to \mathcal{C}$
is a group object if $f([0])$ is a final object in $\mathcal{C}$.
We denote by $\Grp(\mathcal{C})$ the full subcategory of $\Fun(\textup{N}(\Delta)^{op},\mathcal{C})$ that is spanned by group objects in $\mathcal{C}$.
Note that $\Grp(\mathcal{C})$ is a full subcategory of $\textup{Mon}(\mathcal{C})$.
\end{Definition}

\begin{Definition}
A derived group scheme over $R$ is a functor
\[
G:\CAlg_{R}\longrightarrow \Grp(\SSS)
\]
such that the composite $\CAlg_R \to \Grp(\SSS)\to \SSS$
is representable by a derived scheme $X$, where the second map
$\Grp(\SSS)\to \SSS$ is induced by $\{[1]\}\subset\Delta$.
If $X$ is affine, then we shall call it an derived affine group scheme.
\end{Definition}
The $\infty$-category $\Grp(\SSS)$ admits a simple description.
Let $\SSS_{*}$ be the $\infty$-category of pointed spaces.
Namely, $\SSS_{*}$ is the (homotopy) fiber of $\Fun(\Delta^1,\SSS)\to \Fun(\{0\}, \SSS)\simeq \SSS$ over the contractible space $*\in \SSS$.
Let $\SSS_{*,\ge 1}$ be the full subcategory of $\SSS_*$ spanned by
pointed connected spaces.
Then by \cite[7.2.2.11]{HTT} we have a
functor
\[
\SSS_{*,\ge 1}\longrightarrow \Fun(\textup{N}(\Delta)^{op},\SSS_*)
\]
which to any $*\to X\in \SSS_{*,\ge1}$ associates the groupoid
of the \v{C}ech nerve, and it induces an equivalence $\SSS_{*,\ge 1}\simeq \Grp(\SSS_*)$.
Since an initial object in $\SSS_*$ 
is a final object, we easily see that there
is a natural equivalence $\Grp(\SSS_*)\simeq \Grp(\SSS)$
induced by the forgetful functor $\SSS_*\to \SSS$ (cf. \cite[7.2.2.5, 7.2.2.10]{HTT}).
By this identification $\SSS_{*,\ge 1}\simeq \Grp(\SSS)$, the functor
$\Grp(\SSS)\to\SSS$ induced by $[1]\in\Delta$ is equivalent
to the composite
\[
\SSS_{*,\ge 1}\stackrel{\Omega}{\longrightarrow} \SSS_{*}\longrightarrow \SSS
\]
where $\Omega$ is the loop space
functor and the second map is the forgetful functor.
Thus one can say that
a derived group scheme is
a functor $G:\CAlg_{R}\to \SSS_{*,\ge 1}$ such that
the composite \[\CAlg_R \stackrel{G}{\to} \SSS_{*,\ge1}\stackrel{\Omega}{\to}\SSS_*\to \SSS\]
is representable by a derived scheme.

\begin{Remark}
\label{simple}
Note that giving a functor
$G:\CAlg_{R}\to \Fun (\textup{N}(\Delta)^{op},\SSS)$
is equivalent to giving a functor
$G':\textup{N}(\Delta)^{op}\to \Fun(\CAlg_{R},\SSS)$.
Using \cite[5.1.2.3]{HTT} we see that
the condition $G$ factors through $\Grp(\SSS)$ is
equivalent to the condition that
$G'$ is a group object in $\Fun(\CAlg_{R},\SSS)$.
Consequently, we have an equivalence
\[
\Fun(\CAlg_{R},\Grp(\SSS))\simeq\Grp(\Fun(\CAlg_{R},\SSS)).
\]
An object $\Grp(\Fun(\CAlg_{R},\SSS))$ is a derived group scheme if and only if
the image by 
\[
\Grp(\Fun(\CAlg_{R},\SSS))\to \Fun(\CAlg_{R},\SSS)
\]
induced by
$[1]\in \Delta$ is a derived scheme.
Thus a derived group scheme over $R$ is a group object of the $\infty$-category
of derived schemes over $R$.
The $\infty$-category of derived group schemes over $R$
is equivalent to $\Grp(\textup{Sch}_R)$.
\end{Remark}

\subsection{Commutative Hopf ring spectrum}
We focus on the case of derived affine group schemes.
An (usual) affine group schemes is regarded as the Zariski spectrum
of a commutative Hopf-algebra. We will give a similar description
in our setting.
By Remark~\ref{simple}, giving a derived affine scheme is equivalent to
giving a functor $G:\textup{N}(\Delta)\to \CAlg_{R}$ such that
$G^{op}:\textup{N}(\Delta)^{op}\to \CAlg_{R}^{op}=\Aff_R$
is a group object in $\CAlg_{R}^{op}=\Aff_R$.
We regard $G$ as a functor $G:\textup{N}(\Delta)^{op}\to \Aff_R$, which
is a group object.
A monoid object $M:\textup{N}(\Delta)^{op}\to \Aff_R$
is a group object
if and only if
\[
\alpha^* \times \beta^*: M([2]) \to M([1])\times M([1])
\]
is an equivalence where $\alpha:\{0,2\}\hookrightarrow [2]$
and $\beta:\{0,1\}\hookrightarrow [2]$.
We have the natural fully faithful functor
\[
\Grp(\Aff_R)\to \Fun(\NNNN(\Delta),\CAlg_R).
\]
We refer to an object in $\Fun(\NNNN(\Delta),\CAlg_R)$
which lies in the essential image
of this functor as a {\it commutative Hopf ring spectrum} over $R$.
We refer to the essential image, we denote by $\CHopf_R$, as
the $\infty$-category of commutative Hopf ring spectra over $R$.
Note that there is
a natural categorical equivalence
$\CHopf_R^{op}\simeq \Grp(\Aff_R)$, which we refer
to as the $\infty$-category of derived affine group schemes over $R$.
Also, we set $\operatorname{GAff}_R:=\CHopf_R^{op}$.
We refer to an object in the essential image of
$\Fun'(\NNNN(\Delta)^{op},(\CAlg_R)^{op})\subset \Fun(\NNNN(\Delta),\CAlg_R)$
as a commutative bi-ring spectra over $R$.
We remark the standard fact: if $M$ is a monoid object in $\SSS$,
$M$ is a group object in $\SSS$ if and only if
a monoid $\pi_0(M)$ is a group.

\subsection{Derived group schemes, group schemes and examples}
Let $G$ be a derived group scheme over a commutative ring spectrum
$R$.
We will explain how to associate to $G$ a (usual) group scheme $\bar{G}$
over $\pi_0(R)$.
For simplicity, we here assume that $G$ is affine, i.e., $G=\Spec A$.
We impose some conditions on $G$.
Let us suppose either of conditions:

\begin{enumerate}
\renewcommand{\labelenumi}{(\roman{enumi})}

\item $G$ is flat over $R$

\item $A$ and $R$ are connective, that is, $\pi_i(A)=\pi_i(R)=0$ for $i<0$.

\end{enumerate}

We first treat the case (i).
In this case, according to \cite[8.2.2.13]{HA} there is an isomorphism $\pi_0(A\otimes_{R}A)\simeq \pi_0(A)\otimes_{\pi_0(R)}\pi_0(A)$ of commutative rings.
Hence the group object $G:\textup{N}(\Delta)^{op}\to \Aff_R$
induces a group structure
$\bar{G}:\textup{N}(\Delta)^{op}\to \Aff_{\pi_0(R)}^0$ of $\bar{G}=\Spec \pi_0(A)$,
where $\Aff_{\pi_0(R)}^0$ denotes the $\infty$-category of
ordinary affine schemes over $\pi_0(R)$.

Next we consider the case (ii).
In this case, we also have
an isomorphism $\pi_0(A\otimes_{R}A)\simeq \pi_0(A)\otimes_{\pi_0(R)}\pi_0(A)$ of commutative rings.
Thus a similar argument shows that $\bar{G}:=\Spec \pi_0(A)$ inherits a group structure.
In addition,
$\bar{G}$ is equivalent to the composite
\[
G_0:\CAlg_{H\pi_0(R)}^{\textup{dis}} \hookrightarrow \CAlg_R \stackrel{G}{\to} \Grp(\mathcal{S})\stackrel{\pi_0}{\to} \Grp(\mathcal{S}^{\textup{dis}})
\]
where $\CAlg_{H\pi_0(R)}^{\textup{dis}}$ is the nerve of the category
of usual commutative $\pi_0(R)$-rings, the first functor is the natural functor, and $\mathcal{S}^{\textup{dis}}$ is the category of small
sets.
A group scheme $H$ over $\pi_0(R)$ is said to be
the underlying group scheme of a derived group scheme $G$
if $H$ represents the above composite $G_0$.

Conversely, we may regard a flat group scheme $G$ over $\pi_0(R)$
as a derived group scheme that is flat over $H\pi_0(R)$.
Here $H\pi_0(R)$ is the
Eilenberg-MacLane spectrum, which is a discrete commutative
ring spectrum.
Set $G=\Spec B$. Then the usual tensor product
$B\otimes_{\pi_0(R)}B$ of commutative rings coincides with
the ``derived'' tensor product of $HB$ and $HB$ over $H\pi_0(R)$ in
$\CAlg$. Consequently, $G$ can be viewed as
a derived group scheme.
The $\infty$-category of derived affine group schemes over $H\pi_0(R)$
contains the nerve of the category of affine group schemes
which are flat over $\pi_0(R)$ as a full subcategory.

Finally, we give some examples of derived affine group schemes,
which do not necessarily come from usual flat group schemes.

\begin{Example}
\label{barcon}
Let $s:A\to R$ be an augmentation map in $\CAlg_R$.
Then we have a section $s^*:\Spec R\to \Spec A$.
Let $G:=\Spec R\times_{\Spec A}\Spec R$.
The projection morphism $G\times_RG\simeq \Spec R\times_{\Spec A}\Spec R\times_{\Spec A}\Spec R \stackrel{p_{13}}{\longrightarrow} \Spec R\times_{\Spec A}\Spec R\simeq G$ deterrmies a ``multiplication''.
To make this idea precise consider
the \v{C}ech nerve $\NNNN(\Delta)^{op}\to \Aff_R$ associated to $s^*$
(see \cite[6.1.2.11]{HTT}).
This \v{C}ech nerve
is a derived affine group scheme over $R$
whose underlying scheme is $\Spec R\times_{\Spec A}\Spec R$.
In $\CAlg_R$, this construction is known as a bar construction.
\end{Example}

\begin{Example}
Let $R$ be a commutative ring spectrum.
Let $M\in \PMod_R$.
Let $f:\CAlg_R\to \Grp(\SSS)$
be a functor given by $A \mapsto \Aut(M\otimes_RA)$ (see Example~\ref{example1}).
Then according to Lemma~\ref{Homrep}, $f$ is representable by a derived affine
group scheme over $R$. See also Example~\ref{example2}.
\end{Example}

\begin{Example}
Let $\mathbb{S}[\mathbb{CP}^\infty]:=\Sigma^{\infty}\mathbb{CP}^\infty_+$ be the unreduced suspention
spectrum of the classifying space $\mathbb{CP}^\infty$.
The commutative monoid structure in $\mathcal{S}$
(that is, $E_\infty$-structure) of $\mathbb{CP}^\infty$ induces a commutative
algebra structure
on $\mathbb{S}[\mathbb{CP}^\infty]$. Namely, $\mathbb{S}[\mathbb{CP}^\infty]\in \CAlg$.
The diagonal map $\mathbb{CP}^\infty\to \mathbb{CP}^\infty\times \mathbb{CP}^\infty$ makes $\mathbb{S}[\mathbb{CP}^\infty]$ a
commutative Hopf ring spectrum and
thus $\Spec \mathbb{S}[\mathbb{CP}^\infty]$ is a derived affine group scheme over $\mathbb{S}$
(see \cite[12.1]{Ro}).
\end{Example}

\begin{Example}
Let $k$ be a number field.
In \cite{Spi}
Spitzweck constructed the derived affine group scheme $G=\Spec B$
over $H\mathbb{Z}$
such that the $\infty$-category of $H\mathbb{Z}$-spectra with action of $G$
(see below) is equivalent to the stable $\infty$-category of
Voevodsky's category $\mathsf{DM}(k)$ of integer coefficients generated by Tate motives.
(His results are much stronger, see \cite{Spi}.)
\end{Example}

\subsection{$\infty$-categories of commutative bi-ring spectra and commutative
Hopf ring spectra}
We will prove that $\infty$-categories of commutative Hopf ring spectra
and commutative bi-ring spectra have good properties,
that is, these are presentable $\infty$-categories.
To this end, we first give a slighly modified description
of commutative bi-ring spectra.

The $\infty$-category $\CAlg_R$ has the natural coCartesian symmetric
monoidal structure (cf. \cite{HA})
which we will specify by a coCartesian fibration
$\CAlg_R^\otimes \to \NNNN(\FIN)$.
Let $\Ass\to \NNNN(\FIN)$ denote the associative $\infty$-operad
(see \cite[4.1.1.3]{HA} for the definition of associative $\infty$-operad
$\Ass$).
The projection 
\[
p:\CAlg_R^{m\otimes}:=\CAlg_R^\otimes\times_{\NNNN(\FIN)}\Ass\to \Ass
\]
is a monoidal $\infty$-category (cf. \cite[4.1.1.10]{HA}),
that is, the underlying monoidal $\infty$-category of $\CAlg_R^\otimes$.
Let us recall the construction of the opposite monoidal $\infty$-category.
Let $\MMM^\otimes \to \Ass$ be a monoidal $\infty$-category.
Let $F_{\MMM^\otimes}:\Ass\to \wCat$ be a functor corresponding to
$\MMM^\otimes \to \Ass$ via the straightening functor
(see \cite[3.2]{HTT} for the straightening and unstraightening functors).
Let $\operatorname{Op}:\wCat\to \wCat$ be the natural (auto)equivalence which carries $S$ to the opposite
category $S^{op}$.
The composite $\operatorname{Op}\circ F_{\MMM^\otimes}:\Ass\to \wCat$
defines a monoidal $\infty$-category $\MMM_{op}^\otimes \to \Ass$
via the unstraightening functor.
Let $\MMM$ be the underlying $\infty$-category of $\MMM^\otimes$.
Roughly speaking, $\MMM_{op}^\otimes\to \Ass$
is the $\infty$-category $\MMM^{op}$
endowed with the monoidal structure
given by
$\otimes^{op}:(\MMM\times \MMM)^{op}\to \MMM^{op}$ where $\otimes$ indicates the monoidal operation of $\MMM$.
If a monoidal $\infty$-category $\NNN^\otimes\to \Ass$
is equivalent to $\MMM_{op}^\otimes\to \Ass$,
then we shall refer to $\NNN^\otimes\to \Ass$
as the opposite monoidal $\infty$-category of $\MMM^\otimes\to \Ass$.
If we replace $\Ass$ by $\NNNN(\FIN)$, we obtain the opposite
symmetric monoidal $\I$-category of a symmetric monoidal $\I$-category.

Let $q:(\CAlg_R)^{m\otimes}_{op}\to \Ass$ denote the opposite monoidal $\infty$-category of $p$.
Let $\CoAlg(\CAlg_R^\otimes)$ be $\Alg_{/\Ass}((\CAlg_R)^{m\otimes}_{op})^{op}$
where $\Alg_{/\Ass}((\CAlg_R)^{m\otimes}_{op})$ is the $\infty$-category of
algebra objects.
We refer to $\CoAlg(\CAlg_R^\otimes)$ as the $\infty$-category
of commutative bi-ring spectra over $R$ (or commutative bi-ring $R$-module spectra).

Now we show that this definition is compatible with
the above definition.
The opposite symmetric monoidal $\I$-category $(\CAlg_R)^{\otimes}_{op}\to \NNNN(\FIN)$ of $p:\CAlg_R^\otimes\to \NNNN(\FIN)$
is a Cartesian monoidal $\I$-category (cf. \cite[2.4.0.1]{HA}).
By \cite[2.4.1.9]{HA}, there is a Cartesian structure \cite[2.4.1.1]{HA}
$(\CAlg_R)^{\otimes}_{op}\to (\CAlg_R)^{op}$ and it induces the second categorical equivalence in
\[
\Alg_{/\Ass}((\CAlg_R)^{\otimes}_{op})\simeq \Alg_{/\NNNN(\Delta)^{op}}((\CAlg_R)^{\otimes}_{op}\times_{\Ass}\NNNN(\Delta)^{op})\simeq \Fun'(\NNNN(\Delta)^{op},(\CAlg_R)^{op}),
\]
where the first equivalence is induced by the map $\operatorname{Cut}:\NNNN(\Delta)^{op}\to \Ass$ defined in \cite[4.1.2.5]{HA} and \cite[4.1.2.15]{HA},
and $\Fun'(\NNNN(\Delta)^{op},(\CAlg_R)^{op})$
is the full subcategory of monoid objects (see Appendix~\ref{monoidgroup}).
Remark that $f:\NNNN(\Delta)^{op}\to (\CAlg_R)^{op}$ lies in
$\Fun'(\NNNN(\Delta)^{op},(\CAlg_R)^{op})$
if and only if
maps $\{i-1,i\}\hookrightarrow [n]$ with $1\le i\le n$
induce an equivalence
$\otimes_{1\le i\le n}f([1])\to f([n])$ for each $n$, and $f([0])\simeq R$.
Consequently, $\CoAlg(\CAlg_R^\otimes)$ is naturally equivalent to
$\Fun'(\NNNN(\Delta),\CAlg_R)$
where $\Fun'(\NNNN(\Delta),\CAlg_R)$ again denotes
the full subcategory of $\Fun(\NNNN(\Delta),\CAlg_R)$
spanned by comonoid objects.

Let $a=\{0,2\}\hookrightarrow [2]$ and $b=\{0,1\}\hookrightarrow [2]$.
Let $C:\NNNN(\Delta)\to \CAlg_R$ be an object in
$\Fun'(\NNNN(\Delta),\CAlg_R)\simeq \CoAlg(\CAlg_R^\otimes)$.
The object $C$ is a commutative Hopf ring spectrum
if and only if $C(a)$ and $C(b)$ determine $u:C([1])\to C([1])\otimes_RC([1])$
and $v:C([1])\to C([1])\otimes_RC([1])$ such that
$u\otimes v:C([1])\otimes_RC([1]) \to C([1])\otimes_RC([1])$ is an
equivalence in $\CAlg_R$.
The spectrum $R$ is a unit of the symmetric monoidal $\I$-category
$\CAlg_R^{\otimes}$ and thus $R$ is promoted to an object in $\CoAlg(\CAlg_R^\otimes)$. Clearly, $R$ is a commutative Hopf ring spectrum.
The $\infty$-category $\CHopf_R$ is contained in
$\CoAlg(\CAlg_R^\otimes)$ as a full subcategory.
Yoneda lemma implies the natural functor
$\CHopf_R^{op}\to \Fun(\CAlg_R, \Grp(\mathcal{S}))$ is fully faithful.

\begin{Remark}
The natural inclusion $\Fun'(\NNNN(\Delta),\CAlg_R)\to \Fun(\NNNN(\Delta),\CAlg_R)$
preserves small colimits.
Let $I$ be a small $\infty$-category and $I\to \Fun'(\NNNN(\Delta),\CAlg_R)$
a functor.
We will claim that
a colimit of the composition $q:I\to \Fun'(\Delta,\CAlg_R)\to \Fun(\NNNN(\Delta),\CAlg_R)$ satisfies the comonoid condition.
For $\lambda\in I$, we set $A_\lambda=q(\lambda)([1])$.
Note that $q([0])\simeq R$ and $q(\lambda)([n])$ is equivalent to
the $n$-fold tensor product
$A_\lambda\otimes_R\ldots \otimes_RA_\lambda$.
By \cite[5.1.2.3]{HTT}, the $n$-th term of
the colimit of $q$ is $\textup{colim}(A_\lambda\otimes_R\ldots \otimes_RA_\lambda)$ (indexed by $I$) in $\CAlg_R$.
It will suffice to prove
that for each $n\in \mathbb{N}$, inclusions
$\{i-1,i\}\hookrightarrow [n]$ for $1\le i\le n$ induces an equivalence
\[
\textup{colim}(A_\lambda)\otimes_R\ldots \otimes_R\textup{colim}(A_\lambda)\to 
\textup{colim}(A_\lambda\otimes_R\ldots \otimes_RA_\lambda).
\]
According to \cite[4.4.2.7]{HTT}, we may assume that
$I$ is either a pushout diagram or a coproduct diagram.
For simplicity, suppose that $n=2$. (The general case is straightforward.)
Note that the symmetric monoidal structure of $\CAlg_R$ is coCartesian.
In the coproduct case,
$(\otimes_\lambda A_\lambda)\otimes_R(\otimes_\lambda A_\lambda)\simeq \otimes_{\lambda}(A_\lambda\otimes_RA_\lambda)$.
In the case of pushout, for a diagram $A\leftarrow C\to B$ in $\CAlg_R$,
we have an equivalence $(A\otimes_CB)\otimes_R(A\otimes_CB)\simeq (A\otimes_RA)\otimes_{C\otimes_RC}(B\otimes_RB)$.
Hence our claim follows.

The inclusion $\CHopf_R\hookrightarrow \CoAlg(\CAlg_R^\otimes)\simeq \Fun'(\NNNN(\Delta),\CAlg_R)$
preserves small colimits. Let $I$ be a small $\infty$-category
and $I\to \CHopf_R$ a functor.
Let $q:I\to \CHopf_R\hookrightarrow \Fun'(\Delta,\CAlg_R)$ be the composition.
We adopt the notation similar to the above paragraph.
We claim that the colimit of $q$ belongs to $\CHopf_R$.
By assumption, $a=\{0,2\}\hookrightarrow [2]$ and $b=\{0,1\}\hookrightarrow [2]$ and the colimits induce a diagram
\[
\xymatrix{
\textup{colim}A_\lambda \ar[r] \ar[rd]_{\textup{colim}(a_*)} &  \textup{colim}(A_\lambda\otimes_RA_\lambda) \ar[d] & \textup{colim}A_\lambda \ar[l] \ar[ld]^{\textup{colim}(b_*)} \\
   & \textup{colim}(A_\lambda\otimes_RA_\lambda) & 
}
\]
where the the upper horizontal diagram is the colimit of the
coproduct diagrams
$A_\lambda\to A_\lambda\otimes_RA_\lambda\leftarrow A_\lambda$.
The vertical arrow in the middle is an equivalence (by our assumption).
Moreover, in the previous paragraph, we have shown that
the upper horizontal diagram
exhibits $\textup{colim}(A_\lambda\otimes_RA_\lambda)$ as the coproduct
of $\textup{colim}A_\lambda$ and $\textup{colim}A_\lambda$.
This implies that the colimit of $q$ belongs to $\CHopf_R$.
\end{Remark}

\begin{Proposition}
\label{presentability1}
The $\infty$-category $\CoAlg(\CAlg_R^\otimes)$ is 
a presentable $\infty$-category.
\end{Proposition}

\Proof
Let $\mathcal{C}$ be a subcategory of $\wCat$ such that:
\begin{itemize}

\item objects are $\infty$-categories $\mathcal{X}$ such that $\mathcal{X}^{op}$
is an accessible $\infty$-category,

\item morphisms are functors $F:\mathcal{X}\to \mathcal{Y}$
such that $F^{op}:\mathcal{X}^{op}\to \mathcal{Y}^{op}$ are accessible
functors.

\end{itemize}
Note that $\textup{Op}:\wCat\to \wCat$ which sends $\mathcal{X}$ to $\mathcal{X}^{op}$ is a categorical equivalence.
Moreover by \cite[5.4.7.3]{HTT} the limit of accessible $\infty$-categories
in $\wCat$ exists and it is an accessible $\infty$-category.
These observations together with \cite[5.4.4.3, 5.1.2.3]{HTT}
imply that
$\mathcal{C}\subset \wCat$ satisfies the conditions (a), (b), (c)
in \cite[5.4.7.11]{HTT}. Since the monoidal structure on $\CAlg_R$
is compatible with small colimits, combined with \cite[3.2.3.4]{HA}
we can apply \cite[5.4.7.14]{HTT} to deduce that 
$\CoAlg(\CAlg_R^\otimes)$ is 
accessible.
Finally, $\CoAlg(\CAlg_R^\otimes)$ admits small colimits
since $\Fun(\NNNN(\Delta),\CAlg_R)$ is presentable
and  $\CoAlg(\CAlg_R^\otimes)\subset \Fun(\NNNN(\Delta),\CAlg_R)$
is stable under small colimits.
\QED

\begin{Proposition}
\label{presentabilityhopf}
The $\infty$-category $\CHopf_R$ is 
a presentable $\infty$-category.
\end{Proposition}

\Proof
Let $V \rightarrow \NNNN(\Delta)$ denote the inclusion
corresponding to
\[
\xymatrix{
 & a=\{0,2\} \ar[d] \\
b=\{0,1\} \ar[r] & c=[2]
}
\]
where two maps are inclusions.
Namely, $V$ has exactly three objects $a, b, c$,
and non-degenerate maps are $a\to c$ and $b\to c$.
The composition with $V\to \NNNN(\Delta)$
determines a map
$\Fun'(\NNNN(\Delta),\CAlg_R)\to \Fun(V, \CAlg_R)$.
For $p:V\to \CAlg_R$, $p$ induces
$p(a)\otimes p(b)\to p(c)$ since $p(a)\otimes p(b)$
is a coproduct of $p(a)$ and $p(b)$.
By left Kan extension it yields $\Fun(V, \CAlg_R)\to \Fun (\Delta^1,\CAlg_R)$
which carries $p$ to $p(a)\otimes p(b)\to p(c)$,
and we have the composition
$\sigma:\Fun'(\NNNN(\Delta),\CAlg_R)\to \Fun(\Delta^1,\CAlg_R)$.
By the definition of $\CHopf_R$,
we have a homotopy cartesian square
\[
\xymatrix{
\CHopf_R \ar[r] \ar[d] & \Fun'(\NNNN(\Delta),\CAlg_R) \ar[d]^{\sigma} \\
\Fun_{\simeq}(\Delta^1, \CAlg_R) \ar[r]^{\tau} & \Fun(\Delta^1, \CAlg_R)
}
\]
where $\Fun_{\simeq}(\Delta^1, \CAlg_R)$ is the full subcategory
of $\Fun(\Delta^1, \CAlg_R)$ spanned by maps $\Delta^1\to \CAlg_R$
which correspond to equivalences in $\CAlg_R$, and $\tau$
is the inclusion.
Since $\Fun_{\simeq}(\Delta^1, \CAlg_R)\simeq \CAlg_R$,
$\tau$ preserves small colimits.
According to \cite[5.1.2.3]{HTT}, we see that $\sigma$ preserves small 
colimits (by noting
$\Fun'(\NNNN(\Delta),\CAlg_R)\to \Fun(\NNNN(\Delta),\CAlg_R)$
preserves small colimits).
Note that by Proposition~\ref{presentability1}
and \cite[5.4.4.3]{HTT} $\Fun'(\NNNN(\Delta),\CAlg_R)$,
$\Fun(\Delta^1,\CAlg_R)$ and $\Fun_{\simeq}(\Delta^1, \CAlg_R)$ 
are presentable $\infty$-categories (we remark that $\CAlg_R$ is presentable).
Thus by virtue of \cite[5.5.3.13]{HTT} we see that
$\CHopf_R$ is also presentable.
\QED

As a corollary of these results,
we have

\begin{Corollary}
Let $\operatorname{GAff}_R$ be the $\infty$-category of derived
affine group schemes over $R$.
Then $\operatorname{GAff}_R$ has small colimits and limits.
The forgetful functor $\operatorname{GAff}_R\to \Aff_R$
preserves small limits.
\end{Corollary}

\subsection{Representations of commutative bi-ring spectra and commutative Hopf ring spectra}
We will construct a functor
$\CoAlg(\CAlg)\to \wCat$
which carries $B\in \CoAlg(\CAlg)$ to the stable presentable $\I$-category
$\Rep_B$ consisting of spectra
endowed with coaction of $B$.
Informally, $\Rep_B$ is the $\infty$-category of
spectra $N$ endowed with action of the derived monoid
scheme $\Spec B$ which associates
an {\it automorphism} $N\otimes V\stackrel{\sim}{\to} N\otimes V$
to each point $\Spec V\to \Spec B$ with $V\in \CAlg$.
Thus when $B$ does not lie in $\CHopf$, roughly speaking,
$\Mod_B$ (which we are going to define)
does not coincide with the $\infty$-category of
``comodules'' of $B$.
We believe that the notation $\Rep_B$ is little confusing.

Before we define the $\infty$-category $\Rep_B$ for $B\in \CoAlg(\CAlg)$,
we recall the $\infty$-category $\wsCat$ of stable presnetable $\infty$-categories and the functor
$\CAlg\to \wsCat$
which to any $R\in \CAlg$ associates the $\infty$-category $\Mod_R$
of left $R$-module spectra.

Let $\wsCat$ be the $\infty$-category of presentable stable $\infty$-categories
where morphisms are colimit-preserving functors.
(This category is a subcategory of $\wCat$.)
There is a natural symmetric monoidal structure
on $\wsCat$ which commutes with small colimit separately in each variable
(see \cite[II, 4.2]{DAGn} or \cite[6.3.2]{HA}).
For $\mathcal{C}, \mathcal{D}\in\wsCat$, the
tensor product $\mathcal{C}\otimes\mathcal{D}$
has the following universality:
There is a functor $\mathcal{C}\times \mathcal{D}\to \mathcal{C}\otimes \mathcal{D}$ which preserves small colimits separately in each variable,
and
if $\mathcal{E}$ belongs to $\wsCat$ and
$\Fun_c(\mathcal{C}\times \mathcal{D},\mathcal{E})$
denotes the full subcategory of
$\Fun(\mathcal{C}\times \mathcal{D},\mathcal{E})$
spanned by functors which preserve small colimits separately in each variable,
then the composition induces a categorical equivalence
\[
\Fun^{\textup{L}}(\mathcal{C}\otimes\mathcal{D},\mathcal{E})\to \Fun_c(\mathcal{C}\times \mathcal{D},\mathcal{E}),
\]
where $\Fun^{\textup{L}}(-,-)$ on the left side of the equivalence
indicates the full subcategory of $\Fun(-,-)$ spanned by colimit-preserving
functors.
An object $\CAlg(\wsCat)$ can be regarded as a symmetric monoidal
stable presentable $\infty$-category whose tensor product preserves
small colimits separately in each variable.

Let $\LM$ be the $\I$-operad of left modules (see for the
definition \cite[4.2.1.7]{HA}).
Consider the symmetric monoidal $\infty$-category
$\SP^\otimes\to \NNNN(\FIN)$ of spectra.
The natural fibration $\LM\to \Ass$ and its section $\Ass\hookrightarrow \LM$
of $\I$-operads described
in \cite[4.2.1.9, 4.2.1.10]{HA} determine a map
\[
\phi:\textup{LMod}(\SP)=\Alg_{\LM/\Ass}(\SP^\otimes)\to \Alg_{\Ass/\NNNN(\FIN)}(\SP^\otimes).
\]
By \cite[6.3.3.15]{HA} $\phi$ is a
coCartesian fibrarion
(informally
for $R\to R'\in \Alg_{\Ass/\NNNN(\FIN)}(\SP^\otimes)$ and $(R,M)\in \textup{LMod}(\SP)$, $M\to M\otimes_RR'$ is a coCartesian edge lying over it).
Thus the straightening functor
gives rise to $\Alg_{\Ass/\NNNN(\FIN)}(\SP^\otimes)\to \wCat$
which factors through
$\Alg_{\Ass/\NNNN(\FIN)}(\SP^\otimes)\to \wsCat$.
It is extended to a functor between the
$\infty$-categories of commutative algebra objects
\[
\CAlg(\Alg_{\Ass/\NNNN(\FIN)}(\SP^\otimes))\to \CAlg(\wsCat)
\]
(cf. \cite[6.3.5.16]{HA}).
As explained in the proof of \cite[6.3.5.18]{HA},
the unique bifunctor $\Ass \times \operatorname{Comm}^\otimes\to \operatorname{Comm}^\otimes$ of $\infty$-operads
(here the $\infty$-operad $\operatorname{Comm}^\otimes$ is determined by the identity map $\operatorname{Comm}^\otimes:=\NNNN(\FIN)\to \NNNN(\FIN)$)
exhibits $\operatorname{Comm}^\otimes$ as a tensor product of
$\Ass$ and $\operatorname{Comm}^\otimes$. It follows a categorical equivalence
$\CAlg(\SP^\otimes)\to \CAlg(\Alg_{\Ass/\NNNN(\FIN)}(\SP^\otimes))$.
Thus we have
\[
\Theta:\CAlg\longrightarrow \CAlg(\wsCat)
\]
which carries $A$ to $\Mod_R^\otimes$.

Next using $\Theta$, for any $B\in \CoAlg(\CAlg_R^\otimes)$
we will define
an $\infty$-category $\Mod_B$
in a functorial fashion.
Remember that the $\infty$-category
$\CoAlg(\CAlg_R^\otimes)$ is 
equivalent to the $\infty$-category
$\Fun'(\NNNN(\Delta),\CAlg_R)$
of comonoid objects. 
The functor $\Theta$ naturally induces
$\Theta_R:\CAlg_R\simeq \CAlg_{R/}\to \CAlg(\wsCat)_{\Mod_R^\otimes/}$ (the first equivalence follows from \cite[3.4.1.7]{HA}).
Hence composing
$\CoAlg(\CAlg_R^\otimes)\simeq \Fun'(\NNNN(\Delta),\CAlg_{R})$
with it we have 
\[
\CoAlg(\CAlg^\otimes_R)\to \Fun(\NNNN(\Delta),\CAlg(\wsCat)_{\Mod_R^\otimes/}).
\]
Since $\CAlg(\wsCat)_{\Mod^\otimes_R/}$
admits small limits (because it is presentable by \cite[5.5.3.11]{HTT}), there is a right adjoint of
$\CAlg(\wsCat)_{\Mod_R^\otimes/}\to \Fun(\NNNN(\Delta),\CAlg(\wsCat)_{\Mod_R^\otimes/})$ induced by the obvious map
$\NNNN(\Delta)\to \Delta^{0}$.
Namely, the right adjoint
\[
\Fun(\NNNN(\Delta),\CAlg(\wsCat)_{\Mod_R^\otimes/}) \to \CAlg(\wsCat)_{\Mod_R^\otimes/}
\]
sends $\NNNN(\Delta)\to \CAlg(\wsCat)_{\Mod_R^{\otimes}/}$ to its limit.
Combining all together we have
\[
\CoAlg(\CAlg_R^\otimes)\to \Fun(\NNNN(\Delta),\CAlg(\wsCat)_{\Mod_R^\otimes/})\to \CAlg(\wsCat)_{\Mod_R^\otimes/}
\]
and for $B\in \CoAlg(\CAlg_R^\otimes)$ we set its image $\Mod_R^\otimes\to \Mod_B^\otimes\in \CAlg(\wsCat)_{\Mod_R^\otimes/}$ which we refer to
as the $R$-linear symmetric monoidal $\infty$-category of representations
of the commutative
bi-ring spectrum
$B$ (here $R$-linear structure means a symmetric monoidal colimit-preserving
functor $\Mod_R^\otimes \to \Rep_G^\otimes$).
If $G=\Spec B$ is a derived affine group (monoid) scheme over $R$,
then we often write $\Rep_{G}$ for $\Rep_B$.

\begin{Proposition}
\label{bigrep}
The $\infty$-category $\Rep_G$ is a stable presentable $\infty$-category
endowed with a symmetric monoidal structure which preserves small
colimits separately in each variable.
\end{Proposition}

Let $\Rep^\otimes_G$ denote the symmetric monoidal $\infty$-category
of representations of $G$.
The unit $u:\Spec R\to G$ induces a symmetric monoidal functor
$u^*\Rep_G^\otimes\to \Mod_R^\otimes$.
Let $\PRep_G^\otimes$ be the symmetric monoidal full subcategory
of $\Rep_G^\otimes$ spanned by dualizable objects.
An object $M\in \Rep_G$ lies in $\PRep_G$ if and only if
$u^*(M)$ lies in $\PMod_R$.
We refer to $\PRep_G$ as the $\infty$-category
of perfect representations of $G$.
We can easily deduce the following:

\begin{Proposition}
\label{perrep}
The $\infty$-category $\PRep_G$ is a small stable idempotent complete
$\infty$-category
endowed with a symmetric monoidal structure which preserves finite
colimits separately in each variable.
\end{Proposition}

Let $(\CAlg_R)^{op} \hookrightarrow \Fun(\CAlg_R,\widehat{\mathcal{S}})$
be Yoneda embedding, where $\widehat{\mathcal{S}}$ denotes
the $\infty$-category of (not necessarily small) spaces, i.e. Kan complexes.
We shall refer to objects in $\Fun(\CAlg_R,\widehat{\mathcal{S}})$
as presheaves on $\CAlg_R$ or simply functors.
By left Kan extension of $\Theta_R$,
we have a colimit-preserving functor
\[
\overline{\Theta}_R:\Fun(\CAlg_R,\widehat{\mathcal{S}}) \to \CAlg(\wCat)^{op}.
\]
For $X\in \Fun(\CAlg_R,\widehat{\mathcal{S}})$,
we write $\Mod_X^\otimes$ for $\overline{\Theta}_R(X)$.
We denote by $\PMod_X^\otimes$ the full subcategory spanned by
dualizable objects.
Let $G$ be a derived affine group scheme
and let $\psi:\NNNN(\Delta)^{op}\to \Aff_R$
be the corresponding simplicial object.
Let $\NNNN(\Delta)^{op}\stackrel{\psi}{\to} (\CAlg_R)^{op} \hookrightarrow \Fun(\CAlg_R,\widehat{\mathcal{S}})$ be the composition
and let $\mathsf{B}G$ denote the colimit.
Remember $\overline{\Theta}_R(\mathsf{B}G)=\Mod^\otimes_{\mathsf{B}G}\simeq \Rep_G^\otimes$ and $\PMod_{\mathsf{B}G}^\otimes \simeq \PRep_G^\otimes$.

$\bullet$
The author thanks participants of the seminars on HAG at Kyoto university, and SGAD 2011
for helpful conversations related to the subject of this paper.
Also, he thanks the workshop on motives at Tohoku university
in March 2012 where he gave main contents of this paper.
The author is partly supported by Grant-in-aid for Scientific Reseach 23840003,
Japan Society for the promotion of science.

\end{document}